\documentclass[hidelinks]{article}
\usepackage[authoryear,longnamesfirst]{natbib}

\usepackage[utf8]{inputenc}
\usepackage{amsfonts}
\usepackage{amsmath}
\usepackage{amssymb}
\usepackage{graphicx}
\usepackage{xcolor}
\usepackage{caption}
\usepackage[a4paper, total={5in, 9in}]{geometry}
\usepackage{hyperref}
\usepackage{amsthm}

\begin{document}
\title{Approximation bounds for convolutional neural networks in operator learning}

\author{Nicola R. Franco$^{1}$, Stefania Fresca$^{1}$, Andrea Manzoni$^{1}$, Paolo Zunino$^{1}$}
\date{}
\maketitle
\vspace{-1em}
\noindent\textit{
$^{1}$MOX, Department of Mathematics, Politecnico di Milano,\newline
\textcolor{white}{\hspace{0.6em}}Piazza Leonardo da Vinci 32, Milan, 20133, Italy}\\
\newcommand{\fck}[1]{\mathcal{F}_{(#1)}}
\newcommand{\matrixx}[2]{\left[\begin{array}{cc}#1 & #2\\-#2&#1\end{array}\right]}
\newcommand{\zbmatrix}[1]{\matrixx{\text{Re}(#1)}{\text{Im}(#1)}}
\newcommand{\zb}{\mathbf{z}}
\newcommand{\w}{\mathbf{w}}
\newcommand{\tripleleft}{\left[\left[\left[}
\newcommand{\tripleright}{\right]\right]\right]}
\newcommand{\doubleleft}{\left[\left[}
\newcommand{\doubleright}{\right]\right]}
\newcommand{\x}{\mathbf{x}}

\newtheorem{lemma}{Lemma}
\newtheorem{theorem}{Theorem}
\newtheorem{definition}{Definition}
\newcommand{\fcnn}{\mathcal{F}_{\omega}}
\newcommand{\scnn}{\mathcal{S}_{m}}
\newcommand{\I}{\mathbf{i}}
\newcommand{\mup}{\boldsymbol{\mu}}
\newcommand{\proved}{\hfill\square\vspace{1.5em}}
\newcommand{\review}[1]{#1}

\begin{abstract}
\noindent Recently, deep Convolutional Neural Networks (CNNs) have proven to be successful when employed in areas such as reduced order modeling of parametrized PDEs. Despite their accuracy and efficiency, the approaches available in the literature still lack a rigorous justification on their mathematical foundations. Motivated by this fact, in this paper we derive rigorous error bounds for the approximation of nonlinear operators by means of CNN models.
More precisely, we address the case in which an operator maps a finite dimensional input $\mup\in\mathbb{R}^{p}$ onto a functional output $u_{\mup}:[0,1]^{d}\to\mathbb{R}$, and a neural network model is used to approximate a discretized version of the input-to-output map. 
The resulting error estimates provide a clear interpretation of the hyperparameters defining the neural network architecture. All the proofs are constructive, and they ultimately reveal a deep connection between CNNs and the Fourier transform. Finally, we complement the derived error bounds by numerical experiments that illustrate their application.\\
\end{abstract}


\section{Introduction}
Convolutional Neural Networks (CNNs) have become very popular after their tremendous success in computer vision, with applications ranging from image processing to generative models for images generation \citep{dosovitskiy2015learning, sultana}. From a mathematical point of view, image-like data are equivalent to discrete functional signals defined over rectangular domains and vice versa.

Indeed, each continuous function $u:[0,1]^{2}\to\mathbb{R}$ can be discretized in matrix form as
$$\mathbf{U}:=\left[\begin{array}{ccc}
    u(\x_{1,1}) & \dots & u(\x_{1,N}) \\
    \dots & \dots & \dots \\
    u(\x_{N,1}) & \dots & u(\x_{N,N})
\end{array}\right]\in\mathbb{R}^{N\times N},$$
where $\{\x_{i,j}\}_{i,j=1,\dots,N}$ are the vertices of some uniform grid defined over the unit square. In light of this, CNNs have been recently employed for tasks that go beyond computer vision, such as operator learning and reduced order modelling of parameter-dependent PDEs \citep{franco2021deep,fresca2021comprehensive,fresca2021pod,lee2020model,mucke2021reduced}. 

As an example, let $\Omega=(0,1)^{d}$ and assume we are given an operator $\mup\to u_{\mup}$ that maps a finite dimensional input $\mup\in\mathbb{R}^{p}$ onto some functional signal $u_{\mup}:\Omega\to\mathbb{R}$. This is a classical set-up in  parameter-dependent PDE models, where each parameter instance is associated with the corresponding PDE solution. Once a suitable, discrete grid of points  $\{\x_{i_{1},\dots,i_{d}}\}_{i_{1},\dots,i_{d}=1, \dots, N} \subset\Omega$ has been introduced, the operator of interest can be expressed as 
$$\mathbb{R}^{p}\ni\mup\to\mathbf{U}_{\mup}\in\mathbb{R}^{N\times\cdots\times N},$$
denoting by $\mathbf{U}_{\mup}^{i_{1},\dots,i_{d}}\approx u_{\mup}(\x_{i_{1},\dots,i_{d}})$. The final goal is then to construct a neural network model $$\Phi:\mathbb{R}^{p}\to\mathbb{R}^{N\times\cdots\times N},\;\;\text{such that}\;\;
\Phi(\mup)\approx\mathbf{U}_{\mup},$$
with the idea of replacing an operator that is otherwise computationally expensive to evaluate. As previously mentioned, this task can be successfully achieved by CNNs, as they allow to intrinsically account for underlying spatial correlations. 
However, the literature still lacks a comprehensive mathematical analysis and foundation motivating the remarkable performance shown by CNNs, and the role played by each hyperparameter in a CNN model remains unclear.  
In this work, we aim at addressing these critical points, showing rigorous estimates on the error $\mathcal{E}:=\sup_{\mup}\sup_{\mathbf{j}\in\{1,\dots,N\}^{d}} |u_{\mup}(\x_{\mathbf{j}})-\Phi_{\mathbf{j}}(\mup)|$ generated when approximating the operator of interest by means of CNNs.

\subsection*{Literature review}
Neural Networks (NNs) were known to be universal approximators since \cite{cybenko1989approximation}, however the design of effective NN architectures able to preserve desired accuracy properties had not been in-depth investigated until recent years. A substantial improvement was achieved in \cite{yarotsky}, where a rigorous mathematical meaning to structural properties, such as width and depth of a NN model, was first provided. In particular, Yarotsky proved that any $s$-differentiable scalar-valued map $f:[0,1]^{p}\to\mathbb{R}$ can be approximated uniformly with error \review{$0<\varepsilon<1/2$} by some ReLU Deep Neural Network (DNN) with $c\log(1/\varepsilon)$ layers and $c\varepsilon^{-p/s}\log(1/\varepsilon)$ active weights, where $c=c(p,s,f)$ is some constant that depends on the derivatives of $f$. This result was later extended to more general activation functions and different norms, see e.g. \cite{guhring2020error,guhring2021approximation,siegel2022high}, and adapted to the case of CNNs exploiting some algebraic arguments that link dense and convolutional layers, see e.g. \cite{zhou2020universality, he2021approximation}.

However, all these results are limited to the approximation of scalar-valued maps and they are not suited for operator learning. To this end, it is worth to note the following aspect. Assume that we are interested in approximating a vector-valued map $f:[0,1]^{p}\to\mathbb{R}^{n}$, $f(\mup)=[f_{1}(\mup),\dots,f_{n}(\mup)]$, with a DNN model $\Phi$. Clearly, we could exploit the aforementioned results to approximate each $f_{i}$ with some DNN $\phi_{i}$, and then stack together the models to get $\Phi:=[\phi_{1},\dots,\phi_{n}]$. However, with this construction the number of active weights in $\Phi$ would grow linearly with $n$ as $n\to+\infty$. In our context, where we deal with functional outputs and $n=N^{d}$ comes from having discretized $\Omega = (0,1)^{d}$ with a computational grid with $N$ nodes per side, this would be rather undesirable. 
\\\\
Nevertheless, new approaches are now appearing in the literature, in a first attempt to employ NNs for operator learning. Some of these, such as Neural Operators \citep{kovachki2021neural} and DeepONets \citep{lu2021learning}, work with a continuous functional output space, while a second class of approaches relies on a discretization of the output space, see e.g. \cite{kutyniok2021theoretical}. In this work, we focus on the latter family of approaches.

Neural Operators provide a novel framework for building models between infinite-dimensional spaces, and are essentially based on integral operators. Among them, those that have been mostly investigated are Fourier Neural Operators, for which several error estimates have been derived, see, e.g., \cite{kovachki2021universal}. Conversely, DeepONets are a class of models based on a separation of variables approach, which decouples the input parameter and the space variable at output. Error estimates for DeepONets are also available, see \cite{lanthaler2022error}, and they are mostly settled on the aforementioned results in the scalar case and on those by \cite{schwab2019deep} for high-dimensional inputs. 

Besides these methods, deep learning approaches that discretize the functional output space are also available.  
This need usually arises, for instance, when dealing with parameter-dependent PDEs, whose solutions are usually computed through numerical discretization schemes like, e.g., the finite element method. In this case, the functional output space -- usually given by a Sobolev space -- is replaced by a finite-dimensional trial space (e.g., the space of finite elements of degree $r$ built over a triangulation of the spatial domain $\Omega$). Deep learning approaches of this type were proposed in \cite{bhattacharya2020model} and \cite{kutyniok2021theoretical}. The former relies on linear reduction methods to deal with the functional component at output, and it is able to recover mesh independence. The second one, instead, is purely based on DNNs. Both works provide error estimates, most of which are derived by exploiting the results in the scalar-case and projection arguments. \\

This flourishing literature indicates a growing interest aimd at understanding the properties of DNN models and their potential in operator learning. However, at the best of our knowledge, no comprehensive study has yet been proposed for CNN models, despite these latter \review{being} extremely popular in practical applications. One reason might be that CNN architectures can be traced back to sparse versions of dense models, which led researchers to focus on deriving error bounds for DNNs, see e.g. \cite{petersen2020equivalence}. Moreover, CNN models have been mostly studied for handling high-dimensional data at input and not \review{at output}, as in \cite{he2021approximation}. As a consequence, the available literature is left with a missing piece, which is to understand the approximation properties of convolutional layers when reconstructing functional signals. In the present work, we aim at addressing this issue.

\subsection*{Our contribution}
Let $\mup\to u_{\mup}$ be some nonlinear operator whose output are functions $u_{\mup}:[0,1]^{d}\to\mathbb{R}$ defined over the unit hypercube. We provide error bounds for the approximation of such an operator via a CNN model $\Phi:\mathbb{R}^{p}\to\left(\mathbb{R}^{N}\right)^{d}$. In particular, we characterize the model architecture in terms of the approximation error 
$$\mathcal{E}:=\sup_{\mup\in\Theta}\sup_{\mathbf{j}\in\{1,\dots,N\}^{d}} |u_{\mup}(\x_{\mathbf{j}})-\Phi_{\mathbf{j}}(\mup)|,$$
where $\Theta\subset\mathbb{R}^{p}$ is some parameter space and $\{\x_{\mathbf{j}}\}_{\mathbf{j}}\subset[0,1]^{d}$ is a suitable $N\times\dots\times N$ grid. By doing so, we also provide a clear interpretation to the model hyperparameters, including the number of dense and convolutional layers, the amount of active weights and the convolutional channels. In the present work, we limit ourselves to the 1-dimensional case, $d=1$, even if the ideas at the core of our proofs can be extended to higher dimensions with little effort.

We report below our main result, Theorem \ref{th:param_preview}, in which we characterize the approximation error 
in terms of the complexity of a DNN model comprised of a dense and a convolutional block. In what follows, we denote by $H^{s}(\Omega)$ the Sobolev space of $s$-times weakly-differentiable maps with square-integrable derivatives. \review{For a precise mathematical definition about CNN architectures, instead, the reader may refer to the Appendix section at the end of the paper.}
\newtheorem{teorema}{Theorem}
\setcounter{teorema}{1}
\begin{teorema}\label{th:param_preview}
Let $\Omega:=(0,1)$ and let $\{x_{j}\}_{j=1}^{N_{h}}\subset\Omega$ be a uniform grid with stepsize $h=2^{-k}$. Let  $\Theta\subset\mathbb{R}^{p}\ni\mup\to u_{\mup}\in H^{s}(\Omega)$ be  a (nonlinear) operator, where $\Theta$ is a compact domain and $s\ge1$. For some $r\ge0$, assume that the operator is \review{$r$-times continuously Fréchet differentiable. For any $0<\varepsilon<1/2$,} there exists a deep neural network $\Phi:\mathbb{R}^{p}\to\mathbb{R}^{N_{h}}$ such that
$$|u_{\mup}(x_{j})-\Phi_{j}(\mup)|<\varepsilon$$
uniformly for all $\mup\in\Theta$ and all $j=1,\dots,N_{h}$. \review{Additionally, $\Phi$ can be defined as the composition of a fully connected network $\phi$ and a convolutional neural network $\Psi$, i.e. $\Phi=\Psi\circ\phi$, such that the overall architecture has at most}
\begin{itemize}
    \item [i)] $C\log(1/\varepsilon)$ dense layers, with ReLU activation, and $C\log(1/h)$ convolutional layers,
    \item [ii)] $C\varepsilon^{-2/(2s-1)}\left[\review{\varepsilon^{-p/r}}\log(1/\varepsilon)+\log(1/h)\right]$ active\\weights,
    \item [iii)] $C\varepsilon^{-2/(2s-1)}$ channels in input and output,
\end{itemize}
where $C>0$ is some constant dependent on $\Theta$ and on the operator $\mup\to u_{\mup}$, thus also on $s,r,p$.
\end{teorema}

\noindent In particular, the above result shows that:
\begin{enumerate}
    \item[(i)] The number of dense layers depends logarithmically on the desired accuracy, while that of the convolutional layers depends logarithmically on the mesh resolution, i.e. on the number of discretization points.
    \item[(ii)] The width of the dense block is related to the regularity of the operator itself, with smooth operators requiring less neurons.
    \item[(iii)] The number of convolutional features depends on the regularity of the signals $u_{\mup}$ at output.
\end{enumerate}

We mention that, while being the ultimate focus of our work, Theorem \ref{th:param_preview} is only proved at the end of the paper, as we first need to derive some preliminary results. More precisely, the paper is organized as follows. First, in Section \ref{sec:interp} we establish a link between convolutional layers and the Fourier transform. Then, in Section \ref{sec:approx}, we exploit these results to build a CNN model capable of reconstructing any functional output. Finally, in Section \ref{sec:param}, we resort to the parametrized setting and we prove rigorously Theorem \ref{theorem:param}. In addition, we report in Section \ref{sec:valid} some numerical experiments, where we assess the predicted error bounds. A discussion on possible extensions of our results to higher-dimensions can be found in Section \ref{sec:conc}. \review{Finally, in order to keep the paper self-contained, preliminary notions, such as the formal definition of CNN models, are reported in the Appendix.}

\section{Interpolation of Fourier modes}
\label{sec:interp}
Convolution operations are intimately connected to the Fourier transform via the so-called Convolution Theorem, see e.g. \cite{katznelson}. Here, we further investigate this connection by deriving some preliminary results that will serve as building blocks for Theorems \ref{theorem:single} and \ref{theorem:param}. The idea can be stated as follows. Given any dyadic partition of the unit interval, $$\{x_{j}\}_{j=1}^{N_{h}}:=\{(j-1)2^{-k}\}_{j=1}^{N_{h}},$$ where $N_{h}:=2^{k}+1$, and any positive integer $m$, we construct a CNN model $\mathcal{S}_{m}$ that interpolates the (discrete) map
$$[\zb_{-m},\cdots,\zb_{m}]\to\left[\sum_{k=-m}^{m}\zb_{k}\textrm{e}^{2\pi\I k x_{0}},\dots,\sum_{k=-m}^{m}\zb_{k}\textrm{e}^{2\pi\I k x_{N_{h}}}\right]$$
associating the coefficient $\zb_{k}\in\mathbb{C}$, $k = -m, \ldots, m$, to the truncated Fourier transform at the points $x_j$, $j = 1\,, \ldots, N_h$, with $\I$ the imaginary unit. 
The construction of $\mathcal{S}_{m}$ is detailed step-by-step, starting at Lemma \ref{lemma:fundamental} and concluding with Lemma \ref{lemma:fourier}. The proofs are constructive, as they explicitly describe how to implement $\mathcal{S}_{m}$.
In particular, we are able to characterize the complexity of $\mathcal{S}_{m}$ in terms of those specific features that are typical of CNNs, such as depth, kernel size, stride, dilation, padding, and number of input-output channels. For instance, we show in Lemma \ref{lemma:fourier} that the depth of $\mathcal{S}_{m}$ grows logarithmically with the grid resolution, while the active weights grow linearly with $m$. These observations will play a key role for deriving the upper bounds in Theorems \ref{theorem:single} and \ref{theorem:param}.

In what follows, we make use of the embedding $\mathbb{C}\hookrightarrow\mathbb{R}^{4}$, $$\zb\to[\text{Re}(\zb), \text{Im}(\zb), \text{Re}(\zb), \text{Im}(\zb)],$$ to represent complex numbers. This will come in handy when trying to mimic the algebra of elements of $\mathbb{C}$ by using neural networks. With this convention, we also let $\mathbb{C}^{n}\hookrightarrow\mathbb{R}^{4\times n}$ in the obvious way. \review{In particular, when defining CNN architectures of the form $\phi:\mathbb{R}^{4\times n}\to\mathbb{R}^{4\times m}$, we shall write $\phi:\mathbb{C}^{n}\to\mathbb{C}^{m}$. However, this is only a matter of notation: in practice, all the networks considered from now on will never deal with complex values (neither at input or output) but only with the equivalent representation in $\mathbb{R}^{4}$.}

\begin{lemma} For any $k\in\mathbb{N}$ and any $\zb\in\mathbb{C}$ there exists a convolutional neural network $\phi^{k}_{\zb}:\mathbb{C}^{2^{k-1}}\to\mathbb{C}^{2^{k}}$ such that
\begin{itemize}
    \item[i)] it is linear (no activation at any level),
    \item[ii)] it only employs 1D convolutional and reshaping operations,
    \item[iii)] it has an architecture of at most six layers,
    \item[iv)] the input and the output of its convolutional layers have at most 8 channels,
    \item[v)] the kernels of the convolutional layers have size at most equal to 2,
\end{itemize}
and such that
$$\phi^{k}_{\zb}([\w_{1},\dots,\w_{2^{k-1}}]) = [\w_{1},\zb\w_{1},\dots,\w_{2^{k-1}},\zb\w_{2^{k-1}}]$$
for all $\w_{1},\dots,\w_{2^{k-1}}\in\mathbb{C}.$
\label{lemma:fundamental}
\end{lemma}
\noindent\textit{Proof.} Let $n=2^{k-1}$ be the (complex) input dimension. Let $f_{1}$ be a 1D transposed convolutional layer with the following specifics. The layer has four channels at input and four channels at output. It uses a 2-sized window that acts with a stride of 2. The layer has no bias and its weight matrix $\mathbf{W}_{1}\in\mathbb{R}^{4\times 4 \times 2}$, which is obtained by stacking together the convolutional kernels, is zero at all but six entries. These are given by the relations below
$$\left[
\begin{array}{cc}
    \mathbf{W}_{1}^{1,1,1} & \mathbf{W}_{1}^{1,1,2} \vspace{0.5em}\\
    \mathbf{W}_{1}^{2,2,1} & \mathbf{W}_{1}^{2,2,2} \vspace{0.5em}\\
    \mathbf{W}_{1}^{3,3,1} & \mathbf{W}_{1}^{3,3,2} \vspace{0.5em}\\
    \mathbf{W}_{1}^{4,4,1} & \mathbf{W}_{1}^{4,4,2}
\end{array}
\right] = 
\left[
\begin{array}{cc}
    1 & 0 \vspace{0.5em}\\
    0 & 1 \vspace{0.5em}\\
    \text{Re}(\zb) & \text{Im}(\zb) \vspace{0.5em}\\
    -\text{Im}(\zb) & \text{Re}(\zb) 
\end{array}
\right].
$$
Note that above we are also listing some of the zero entries in $\mathbf{W}_{1}$. In this way, it is easier to see what the purpose of $f_{1}$ is. The first block in $\mathbf{W}_{1}$ is used to mimic the action of the identity matrix. Conversely, the second block encodes a $2\times 2$ matrix representation of the complex number $\zb$. The idea is that these two blocks should provide a way of computing the map $\w\to[\w,\zb\w]$. However, for these computations to be actually carried out, we also need a further layer that performs a suitable summation of the outputs given by $f_{1}$.
To this end, we define the second layer, $f_{2}$, as a 1D convolution that maps 4-channeled inputs onto 2-channeled outputs. The latter uses a 1-sized window that acts with a stride of 1. The layer has no bias and its weight matrix $\mathbf{W}_{1}\in\mathbb{R}^{2\times 4 \times 1}$ contains either zeros or ones. The positive entries are
$$\mathbf{W}_{2}^{1,1,1},\;\mathbf{W}_{2}^{1,2,1},\;\mathbf{W}_{2}^{2,3,1},\;\mathbf{W}_{2}^{2,4,1}=1.$$
Then, $f_{2}\circ f_{1}:\mathbb{R}^{4\times n}\to \mathbb{R}^{2 \times 2n}$, and, upto some basic calculations, we have
$$(f_{2}\circ f_{1})\left(
\left[
\w_{1},
\dots,
\w_{n}
\right]\right)=$$
$$\left[
\begin{array}{lllll}
\text{Re}(\w_{1}), & \text{Im}(\w_{1}), & \dots, &  \text{Re}(\w_{n}), & \text{Im}(\w_{n})\\
\text{Re}(\zb\w_{1}), & \text{Im}(\zb\w_{1}), & \dots, &  \text{Re}(\zb\w_{n}), & \text{Im}(\zb\w_{n})
\end{array}
\right].
$$
In practice, the desired output of $\phi_{\zb}^{k}$ is already there, but we need to adjust the output dimension in order to match our convention for complex numbers. 

To this end, we start by introducing a reshape operation $R_{1}:\mathbb{R}^{2\times 2n}\to\mathbb{R}^{1 \times 4n}$ that flattens the whole output. Then, we add a third convolutional layer, $f_{3}$, whose purpose is to double the entries in input. More precisely, we define $f_{3}$ has a 1D convolution that has 1 channel at input and 4 at output. The layer uses a 2-sized kernel that acts with a stride of 2. Once again, the layer introduces no bias and has a weight matrix $\mathbf{W}_{3}\in\mathbb{R}^{4\times1\times2}$ given by
$$\mathbf{W}_{3}=[[[1,0]],\;[[0,1]],\;[[1,0]],\;[[0,1]]].$$
With the notation adopted to represent complex numbers, the current action of $f_{3}\circ R_{1}\circ f_{2}\circ f_{1}$ becomes
$$
\left[
\w_{1},
\dots,
\w_{n}
\right]\to
\left[
\w_{1}, \dots, \w_{n}, \zb \w_{1},\dots, \zb \w_{n} 
\right].
$$
Let us now act further on the output to sort the entries in the desired order. To do so, we introduce a 1D convolutional layer, $f_{4}$, that has a dilation factor of $2^{k}$ (this is because we want to group $\w_{1}$ with $\zb \w_{1}$, which is $2^{k}$ entries faraway, and so on). We let $f_{4}$ go from 4 to 8 channels, and employ a kernel of size 2 with unit stride. Once again, $f_{4}$ does not have a bias term, while its weight matrix
satisfies
$$\mathbf{W}_{4}^{i,j,k} = \begin{cases} 1  & \text{if }i = j+4(k-1)\\0 &\text{otherwise}.\end{cases}$$
At this point we have,
$$f_{4}\circ f_{3}\circ R_{1}\circ f_{2}\circ f_{1}: 
\left[
\w_{1},
\dots,
\w_{n}
\right]\to
\left[
\begin{array}{lll}
 \w_{1},       & \dots,  & \w_{n}\\
 \zb\w_{1},    & \dots,  & \zb\w_{n}
\end{array}
\right],
$$
and we only need to add a final reshaping operation $R_{2}$. We let $R_{2}$ to act as follows. First, it transposes the input by mapping $\mathbb{R}^{8\times n}\to\mathbb{R}^{n\times 8}$. Then, it performs the reshaping $\mathbb{R}^{n\times 8}\to\mathbb{R}^{2n\times 4}$, where entries are read by rows, and finally it transposes back the input so that it ends up in $\mathbb{R}^{4\times 2n}\cong\mathbb{C}^{2^{k}}$. Finally, letting $\phi_{\zb}^{k}:=R_{2}\circ f_{4}\circ f_{3}\circ R_{1}\circ f_{2} \circ f_{1}$ concludes the proof.
$\proved$

\begin{lemma}\label{lemma:modes} Let $k\in\mathbb{N}$ and $h=2^{-k}$. Let $\{x_{1},\dots,x_{N_{h}}\}$ be an uniform partition of $[0,1]$ with stepsize $h$. For any $\omega\in\mathbb{R}$, there exists a convolutional neural network $\fcnn:\mathbb{C}\to\mathbb{C}^{N_{h}-1}$ such that
\begin{itemize}
    \item[i)] $\fcnn$ is linear (no activation at any level),
    \item[ii)] $\fcnn$ has at most depth $C\log(1/h),$
    \item[iii)] $\fcnn$ has at most $C\log(1/h)$ active weights,
    \item[iv)] $\fcnn(\w)=[\w\textrm{e}^{\I\omega x_{1}},\dots,\w\textrm{e}^{\I\omega x_{N_{h}-1}}]$.
\end{itemize} 
\noindent where $C>0$ is a constant independent on $h$ and $\omega$. Furthermore, up to reshape operations, $\fcnn$ only uses 1D convolutional layers that have at most 8 channels at both input and output. Moreover, the kernel size of all convolutional layers in $\fcnn$ is at most 2.\end{lemma}
\noindent\textit{Proof}. Let $\zb\in\mathbb{C}$. For all $j=1,\dots, k$, define the CNNs $\phi^{j}_{\zb^{2^{j-1}}}$ as in Lemma 1. For the sake of simplicity assume that $k\ge2$. Then it is straightforward to check that $\phi_{\zb}^{1}:\w\to[\w,\w\zb]$, while $\phi_{\zb^{2}}^{2}\circ\phi_{\zb}^{1}:\w\to[\w,\w\zb^{2},\w\zb,\w\zb^{3}]$
and so on. In particular, 
$$\phi_{\zb^{2^{k-1}}}^{k}\circ\dots\circ\phi_{\zb}^{1}:\w\to\review{\pi_{k}}\left([\w,\w\zb,\w\zb^{2},\dots,\w\zb^{2^{k}}]\right)$$
where $\review{\pi_{k}}:\mathbb{C}^{N_{h}-1}\to\mathbb{C}^{N_{h}-1}$ is some (invertible) map that acts as a permutation over the entries. \review{We now claim that this permutation can be nullified through the composition of three suitable reshape operations: that is, there exist three reshape layers $R_{1}^{(k)},R_{2}^{(k)},R_{3}^{(k)}$ such that \begin{equation}
    \label{eq:reshape}
    \tag{$\sharp$}
    \pi_{k}^{-1}=R_{3}^{(k)}\circ R_{2}^{(k)}\circ R_{1}^{(k)}.
\end{equation}Before proving it, we note that \eqref{eq:reshape} would immediately yield the desired conclusion. In fact, if we let $\zb:=\textrm{e}^{\I\omega h}$, then \eqref{eq:reshape} allows us to define $\mathcal{F}_{\omega}:=\pi_{k}^{-1}\circ\phi_{\zb^{2^{k-1}}}^{k}\circ\dots\circ\phi_{\zb}^{1}$, which results in the map
$$\mathcal{F}_{\omega}:\w\to[\w\textrm{e}^{\I\omega 0h},\w\textrm{e}^{\I\omega 1h},\dots,\w\textrm{e}^{\I\omega (N_{h}-2)h}].$$
Therefore, it is sufficient for us to prove that \eqref{eq:reshape} holds. To this end, we define the three maps as
$$R_{1}^{(k)}:\mathbb{C}^{2^{k}}\to\mathbb{C}^{\footnotesize{\overbrace{2\times\dots\times2}^{k\textnormal{ times}}}},$$
where the output entries are filled by rows,
$$R_{2}^{(k)}:\mathbb{C}^{2\times\dots\times2}\to\mathbb{C}^{2\times\dots\times2},$$
which acts as a full transposition of the indexes, that is
$$R_{2}^{(k)}(\mathbf{Q}_{i_{1},\dots,i_{k}})=\mathbf{Q}_{i_{k},\dots,i_{1}},$$
and finally $$R_{3}^{(k)}:\mathbb{C}^{2\times\dots\times2}\to\mathbb{C}^{2^{k}},$$ that flattens back the input. With this setup, we shall now prove \eqref{eq:reshape} following an argument by induction over $k$. To start, let $k=2$. We choose to skip the trivial case $k=1$ to provide the reader with a more insightful computation that anticipates the ideas used later for the inductive step. In this case, we note that  $\pi_{k}([\mathbf{0},\mathbf{1},\mathbf{2},\mathbf{3}])=[\mathbf{0},\mathbf{2},\mathbf{1},\mathbf{3}]$, where we write $\mathbf{k}=k+0\I\in\mathbb{C}$ to embed real numbers in $\mathbb{C}$. The reshape layers act on the latter vector as
\begin{multline*}
R_{3}^{(k)}\circ R_{2}^{(k)}\circ R_{1}^{(k)}\left([\mathbf{0},\mathbf{2},\mathbf{1},\mathbf{3}]\right)=
\\
=
R_{3}^{(k)}\circ R_{2}^{(k)}\left(\left[
\begin{array}{cc}
    \mathbf{0} & \mathbf{2} \\
    \mathbf{1} & \mathbf{3}
\end{array}\right]\right)=
R_{3}^{(k)}\left(\left[
\begin{array}{cc}
    \mathbf{0} & \mathbf{1} \\
    \mathbf{2} & \mathbf{3}
\end{array}\right]\right)=\\=
[\mathbf{0},\mathbf{1},\mathbf{2},\mathbf{3}].
\end{multline*}
Since $\pi_{k}$ is a permutation, the above directly implies \eqref{eq:reshape}. To conclude, we shall now prove the inductive step: assuming that \eqref{eq:reshape} holds for $k$, we show that the statement is also true for $k+1$. For any $j\in\mathbb{N}$, let
$$\mathbf{E}_{j}:=\pi_{j}\left(\left[\mathbf{0},\dots,\mathbf{2}^{j}-\mathbf{1}\right]\right).$$ We split the vector $\mathbf{e}_{k+1}$ into halves by introducing the following notation
$$\mathbf{E}_{k+1}=[\mathbf{A}_{1},\dots,\mathbf{A}_{2^{k}},\mathbf{B}_{1},\dots,\mathbf{B}_{2^{k}}],$$
so that $\mathbf{E}_{k+1}$ is the concatenation of the two (complex) vectors $\mathbf{A},\mathbf{B}\in\mathbb{C}^{2^{k}}$. By the construction of Lemma \ref{lemma:fundamental}, it is straightforward to see that
\begin{equation}
    \label{eq:split}
    \tag{$\sharp\sharp$}
    \mathbf{A}=2\mathbf{E}_{k},\quad\mathbf{B}=2\mathbf{E}_{k}+1,
\end{equation}
where, with little abuse of notation, we intend $[\mathbf{v}_{1},\dots,\mathbf{v}_{n}]+1:=[\mathbf{v}_{1}+1,\dots,\mathbf{v}_{n}+1].$
Let now $$\mathbf{Z}:=R_{1}^{(k+1)}(\mathbf{e}_{k+1}),\quad\quad\tilde{\mathbf{Z}}:=R_{2}^{k+1}(\mathbf{Z}).$$ 
By definition, we have
$$\tilde{\mathbf{Z}}_{j_{1},\dots,j_{k},1}=\mathbf{Z}_{1,j_{k},\dots,j_{1}}=\mathbf{A}_{(j_{k}-1)2^{1-1}+\dots+(j_{1}-1)2^{k-1}},$$
$$\tilde{\mathbf{Z}}_{j_{1},\dots,j_{k},2}=\mathbf{Z}_{2,j_{k},\dots,j_{1}}=\mathbf{B}_{(j_{k}-1)2^{1-1}+\dots+(j_{1}-1)2^{k-1}}.$$
We note that, the subtensor of $\tilde{\mathbf{Z}}$ obtained by fixing the last index equal to 1 is equivalent to $R_{2}^{(k)}\circ R_{1}^{(k)}(\mathbf{A})$, and similarly for $\mathbf{B}$. In particular, if we let $$\hat{\mathbf{A}}:=R_{3}^{(k)}\circ R_{2}^{(k)}\circ R_{1}^{(k)}(\mathbf{A}),\quad\hat{\mathbf{B}}:=R_{3}^{(k)}\circ R_{2}^{(k)}\circ R_{1}^{(k)}(\mathbf{B}),$$
then the action of the final layer, $R_{3}^{(k+1)}$, results in
$$
R_{3}^{(k+1)}(\tilde{\mathbf{Z}})=[\mathbf{\hat{A}}_{1},\mathbf{\hat{B}}_{1},\dots,\mathbf{\hat{A}}_{2^{k}},\mathbf{\hat{B}}_{2^{k}}].
$$
Finally, by recalling \eqref{eq:split} and by applying our inductive hypothesis, we have
$$\hat{\mathbf{A}}=[\mathbf{0},\mathbf{2},\dots,\mathbf{2^{k+1}}-\mathbf{2}],\quad\hat{\mathbf{B}}=[\mathbf{1},\mathbf{3},\dots,\mathbf{2^{k+1}}-\mathbf{1}]$$
$$\implies R_{3}^{(k+1)}(\tilde{\mathbf{Z}}) = [\mathbf{0},\mathbf{1},\dots,\mathbf{2^{k+1}}-\mathbf{1}]=\pi_{k+1}^{-1}(\mathbf{E}_{k+1}),$$
which proves our original claim in \eqref{eq:reshape}.
}
$\proved$

\begin{lemma}\label{lemma:fourier} Let $k\in\mathbb{N}$ and $h=2^{-k}$. Let $\{x_{1},\dots,x_{N_{h}}\}$ be a uniform partition of $[0,1]$ with stepsize $h$. For any positive integer $m$, there exists a convolutional neural network $\scnn:\mathbb{C}^{2m+1}\to\mathbb{C}^{N_{h}}$ such that
\begin{itemize}
    \item[i)] $\scnn$ is linear (no activation at any level),
    \item[ii)] $\scnn$ has at most depth $C\log(1/h),$
    \item[iii)] $\scnn$ has at most $Cm\log(1/h)$ active weights,
    \item[iv)] $\scnn$ uses convolutional layers with at most $Cm$ channels,
    \item[v)] for any complex vector $\mathbf{Z}=[\zb_{-m},\dots,\zb_{m}]\in\mathbb{C}^{2m+1}$ one has
    $$\scnn(\mathbf{Z})_{i}=\sum_{k=-m}^{m}\zb_{k}\textrm{e}^{2\pi\I k x_{i}},$$
    for all $i=1,\dots, N_{h}$, where $\scnn(\mathbf{Z})_{i}$ is the $i$th component of the output vector $\scnn(\mathbf{Z})$.
\end{itemize} 
\noindent Here, $C>0$ is a universal constant independent on $h$ and $m$. Furthermore, up to reshape operations, $\scnn$ only uses 1D convolutional layers whose kernel size does not \review{exceed 2}.\end{lemma}
\noindent\textit{Proof}. Fix any $k\in\{-m,\dots,m\}.$ Let $\fck{k}$ be the CNN in Lemma \ref{lemma:modes} when $\omega=2\pi k$. We note that, as $k$ varies, the structure of $\fck{-m},\dots\fck{m}$ does not change: these architectures have the same depth and they employ convolutional layers with the same specifics. Also, the reshaping operations entailed by the networks occur at the same locations. 

Therefore, we can stack all these models on top of each other to obtain a global CNN $\tilde{\mathcal{S}}_{m}$ such that
$$\tilde{\mathcal{S}}_{m}(\mathbf{Z})=\left[\fck{-m}(\zb_{-m}),\dots,\fck{m}(\zb_{m})\right]$$
where $\mathbf{Z}=[\zb_{-m},\dots,\zb_{m}]\in\mathbb{C}^{2m+1}$ is a generic input vector. This can be done as follows. 

To stack $2m+1$ convolutional layers with $c_{in}$ channels at input and $c_{out}$ channels at output each, we define a single CNN layer with $(2m+1) c_{in}$ channels at input and $(2m+1) c_{out}$ at output. Then, to avoid the introduction of redundant kernels, we constrain the new layer to group its kernels in subsets of $(2m+1)$. This ensures the wished behavior, i.e. that we actually stack the outputs of the $2m+1$ original layers as if they work in parallel (thus each seeing only the part of interest of the input). Similarly, reshaping and transpositions can be easily stacked together. For instance, stacking $(2m+1)$ transpositions of the form $\phi:\mathbb{R}^{a\times b}\to\mathbb{R}^{b\times a}$ results in a map from $\mathbb{R}^{(2m+1)\times a \times b}$ to $\mathbb{R}^{(2m+1)\times b\times a}.$

Since $\tilde{\mathcal{S}}_{m}$ takes values in $\mathbb{C}^{(2m+1)\times(N_{h}-1)}$, our next purpose is to append a further layer $L$ such that
$$(L\circ\tilde{\mathcal{S}}_{m})(\mathbf{Z})=\sum_{k=-m}^{m}\fck{k}(\zb_{k})\in\mathbb{C}^{N_{h}-1}.$$ 
It is easy to see that $L$ can be obtained with a convolutional layer having $(2m+1)$ channels at input and $1$ at output, no stride or dilation, and a kernel of size 1 whose weight is constantly equal to 1. 

Finally, we note that for all $k\in\{-m,\dots,m\}$ we have 
$$\textrm{e}^{2\pi\I k x_{1}}=\textrm{e}^{2\pi\I k 0}=1=\textrm{e}^{2\pi\I k}=\textrm{e}^{2\pi\I k x_{N_{h}}}$$
due to periodicity. Therefore, we may simply define $\mathcal{S}_{m}:= A\circ L\circ \tilde{S}_{m}$, where $A$ has the only purpose of appending a copy of its first output at the end, that is
$$A(\w_{1},\dots,\w_{N_{h}-1})=[\w_{1},\dots,\w_{N_{h}-1},\w_{1}].$$ This can be seen as a form of reshaping or padding.

By construction, $\mathcal{S}_{m}$ satisfies (v). Similarly, (i) and (iv) hold. In fact, each of the $\fck{k}$ has length $C\log(1/h)$, where $C$ is a common constant. Since we stacked them in parallel to get $\tilde{\mathcal{S}}_{m}$, our final model has depth $C\log(1/h)+2=\tilde{C}\log(1/h)$. Also, the CNNs $\fck{k}$ featured at most 8 channels, thus $\scnn$ uses no more than $(2m+1)8=\tilde{C}m$ channels at input-output. Property (iii) follows similarly by recalling that we grouped the CNNs kernels in order to properly stack the architectures.$\proved$

\section{Approximation of a single function}
\label{sec:approx}
We shall now exploit the theory developed in Section \ref{sec:interp} in order to derive suitable error bounds for CNNs. To this end, our primal goal is to show that, for any desired accuracy, there exists a single CNN architecture that is able to provide as output an approximation for any function belonging to a given class in terms of smoothness. More precisely, let $s\ge1$ be a smoothness index, and fix some $m\in\mathbb{N}$. Let also $\{x_{j}\}_{j=1}^{N_{h}}$ be some dyadic grid defined over $\Omega=(0,1)$. We build a CNN model $\Psi:\mathbb{C}^{2m+1}\to\mathbb{R}^{N_{h}}$ such that
$$\forall f\in H^{s}(\Omega)\;\;\exists\mathbf{Z}_{f}\in\mathbb{C}^{2m+1}\;\;\text{such that}$$$$\sup_{j=1,\dots,N_{h}}|f(x_{j})-\Psi_{j}(\mathbf{Z}_{f})|< Cm^{1/2-s}||f||_{H}^{s}(\Omega),$$
where $C=C(s)>0$ is some constant and $\Psi_{j}$ is the $j$th component of the CNN output. The above states that any smooth function $f$ can be well approximated by $\Psi$, provided that the model is fed with a suitable input vector. As for Lemma \ref{lemma:fourier}, we characterize the network complexity in terms of depth, channels and active weights. Furthermore, we show that the map $f\to\mathbf{Z}_{f}$ can be realized by some continuous linear operator that depends, at most, on $s$. Before stating this rigorously in Theorem \ref{theorem:single}, it is worth to remark that this result concerns the approximation of any functional output in $H^{s}(\Omega)$. In particular, although the proof is based on classical estimates coming from the literature of Fourier series, no periodicity is required.

\begin{theorem}
\label{theorem:single}
Let $k\in\mathbb{N}$ and $h=2^{-k}$. Let $\{x_{j}\}_{j=1}^{N_{h}}$ be a uniform partition of $\Omega:=(0,1)$ with stepsize $h$, so that $N_{h}=2^{k}+1$. For any positive integer $m$ and a universal constant $C$ independent on $m$ and $h$, there exists a linear convolutional neural network $\Psi:\mathbb{C}^{2m+1}\to\mathbb{R}^{N_{h}}$ with
\begin{itemize}
    \item[i)] at most $C\log(1/h)$ layers,
    \item[ii)] at most $Cm\log(1/h)$ active weights,
    \item[iii)] at most $8m$ channels in input and output, with kernels grouped by $m$,
\end{itemize}

\noindent such that for any $s\ge1$ and all $f\in H^{s}(\Omega)$ one has
$$\sup_{j=1,\dots,N_{h}}\;|f(x_{j})-\Psi_{j}(Tf)|\le cm^{1/2-s}||f||_{H^{s}(\Omega)}.$$
Here, $c=c(s)>0$ and $T:H^{s}(\Omega)\to\mathbb{C}^{2m+1}$ are  a positive constant and a continuous linear operator that depend on $s$, respectively.
\end{theorem}
\newcommand{\lastchanges}[1]{#1}
\noindent\textit{Proof}. Let us denote by $\mathbb{T}$ the 1-dimensional torus, $\mathbb{T}:=\mathbb{R}/\mathbb{Z}$, so that the spaces $\mathcal{C}^{k}(\mathbb{T})$ refer to those functions that are $k$-times differentiable on the torus, namely
$$\mathcal{C}^{k}(\mathbb{T})=\left\{f\in\mathcal{C}^{k}[0,1],\;f^{(j)}(0)=f^{(j)}(1)\;\;\forall j=1,\dots,k\right\}.$$
We start by defining an operator $T_{0}: H^{s}(\Omega)\to H^{s}(\Omega)$ \lastchanges{that perturbs the signal at input to match suitable boundary conditions.} To this end, we recall that there exist polynomials \lastchanges{$p_{0},\dots,p_{s-1}$ and $q_{0},\dots,q_{s-1}$, of degree $2s-1$,} such that
$$p_{j}^{(k)}(0)=\delta_{j,k},\quad p_{j}^{(k)}(1)=0,$$
$$\lastchanges{q_{j}^{(k)}(0)=0,\quad q_{j}^{(k)}(1)=\delta_{j,k}},$$
for $\lastchanges{j,k\in\{0,\dots,s-1\}},$ see e.g. \cite{hermite}. Then, the linear operator
$$\lastchanges{\mathbf{v}\to P\mathbf{v}:=\sum_{j=0}^{s-1}v_{j} (p_{j}- q_{j})}$$
maps any input vector $\mathbf{v}\in\mathbb{R}^{s-1}$ into a \lastchanges{smooth polynomial with given boundary values}.
By recalling that $H^{s}(\Omega)\hookrightarrow\mathcal{C}^{s-1}(\Omega)$ thanks to classical Sobolev inequalities, we are then allowed to define
$$T_{0}:f\to f + P[f^{(0)}(1)-f^{(0)}(0),\dots,f^{(s-1)}(1)-f^{(s-1)}(0)]$$
so that $\lastchanges{T_{0}: H^{s}(\Omega)\to H^{s}(\Omega)}$. We now introduce the following notation \lastchanges{for "periodicized signals"}. For any $f\in H^{s}(\Omega)$ we let $\tilde{f}$ be defined as
\begin{equation}
    \tilde{f}(x) = 
    \begin{cases}
    (T_{0}f)(2x) & 0\le x \le 1/2\\
    f(2x-1) & 1/2 < x \le 1
    \end{cases}.
\end{equation}
It is straightforward to see that $\tilde{f}\in H^{s}(\Omega)\cap\mathcal{C}^{s-1}(\mathbb{T})$. \lastchanges{For instance,
\begin{multline*}
\tilde{f}(0)=f(0)+\sum_{j=0}^{s-1}[f^{(j)}(1)-f^{(j)}(0)]\cdot[p_{j}(0)-q_{j}(0)]\\=f(0)+[f(1)-f(0)]p_{0}(0)=f(1)=\tilde{f}(1),
\end{multline*}
while
\begin{multline*}
\tilde{f}(1/2)=f(1)+\sum_{j=0}^{s-1}[f^{(j)}(1)-f^{(j)}(0)]\cdot[p_{j}(1)-q_{j}(1)]\\=f(1)-[f(1)-f(0)]q_{0}(1)=f(0)=\lim_{x\to\frac{1}{2}^{+}}\tilde{f}(x).
\end{multline*}}
\;\\\lastchanges{Similar calculations hold for the derivatives as well.
Furthermore,} the mapping $f\to\tilde{f}$ is linear and continuous, in the sense that for some constant $C>0$, depending only on $s$, we have
\begin{equation}
    \tag{*}
    \|\tilde{f}\|_{H^{s}(\Omega)}\le C\|f\|_{H^{s}(\Omega)}
\end{equation}
for all $f\in H^{s}(\Omega)$. \lastchanges{With this construction}, for any positive integer $m$, let $S_{m}\tilde{f}$ be the $m$ truncated Fourier series of the function $\tilde{f}$,
$$\left(S_{m}\tilde{f}\right)(x)=\sum_{j=-m}^{m}c_{\tilde{f}}^{k}e^{2\pi\I k x}$$
where
$$c_{\tilde{f}}^{k}:=\int_{\Omega}\tilde{f}(x)e^{-2\pi\I k x}dx.$$
Since $\tilde{f}\in\mathcal{C}^{s-1}(\mathbb{T})$ and its $s$-derivative is in $L^{2}(\mathbb{T})$, by exploiting classical estimates of Fourier analysis, we have the error bound
\begin{equation}
    \tag{**}
    \|\tilde{f}- S_{m}\tilde{f}\|_{L^{\infty}(\Omega)}\le\sqrt{\frac{2}{2s-1}}m^{1/2-s}\|\tilde{f}\|_{H^{s}(\Omega)}.   
\end{equation}
Let now $T:H^{s}(\Omega)\to\mathbb{C}^{2m+1}$ be defined as
$$T: f\to \left[c_{\tilde{f}}^{-m},\dots,c_{\tilde{f}}^{m}\right],$$
so that $T$ maps each signal into the Fourier coefficients of its periodic alias. Let $\{y_{0},\dots,y_{2N_{h}-1}\}$ be a uniform partition of (0,1) that is twice as fine as the original one $\{x_{0},\dots, x_{N_{h}}\}$, that is $y_{j+1}-y_{j}=h/2$. With this partition as a reference, let then $\mathcal{S}_{m}$ be the CNN in Lemma \ref{lemma:fourier}. By definition, we have
$$\left(S_{m}\tilde{f}\right)(y_{i})=\mathcal{S}_{m}(Tf)_{i}.$$
Finally, let $\Psi:= E\circ R\circ\mathcal{S}_{m}$, where
\begin{itemize}
    \item $R$ is a reshape truncation layer,
    $$R(\w_{1},\dots,\w_{2N_{h}-1}) = [\w_{N_{h}},\dots,\w_{2N_{h}-1}],$$
    that we use to remove the undesired output. Note in fact that, the signal $\tilde{f}$ over $(1/2,1)$ is practically $f$ over $(0,1)$. Thus, in light of (**), we are only interested in the second half of the output.
    \item $E:\mathbb{C}^{N_{h}}\to\mathbb{R}^{N_{h}}$ is the embedding that only keeps the real part of the input. This can also be seen as a reshape layer with a truncation at the end. Since $f$, and thus $\tilde{f}$, are real valued, so are $S_{m}\tilde{f}$ and $\mathcal{S}_{m}(Tf)$. Therefore, we are not losing any information.
\end{itemize}

\noindent Finally, let $j\in\{1,\dots,N_{h}\}$. We have 
\begin{multline*}
    |f(x_{j}) - \Psi_{j}(Tf)| = |f(x_{j}) - \mathcal{S}_{m}(Tf)_{2j-1}| \\ =  |\tilde{f}(y_{2j-1}) - \mathcal{S}_{m}(Tf)_{2j-1}| \\= \left|\tilde{f}(y_{2j-1}) - \left(S_{m}\tilde{f}\right)(y_{2j-1})\right|.
\end{multline*}
Thus, by putting together (**) and (*) we get
\begin{multline*}|f(x_{j}) - \Psi_{j}(Tf)|\le \dots \le C m^{1/2-s}\|\tilde{f}\|_{H^{s}(\Omega)}\\\le C m^{1/2-s}\|f\|_{H^{s}(\Omega)}.\end{multline*}
$\proved$
\renewcommand*{\thefootnote}{$\dagger$}

We remark that, as for the results in the previous Section, the proof is constructive.
\begin{figure*}
    \centering
    \includegraphics[width=\textwidth]{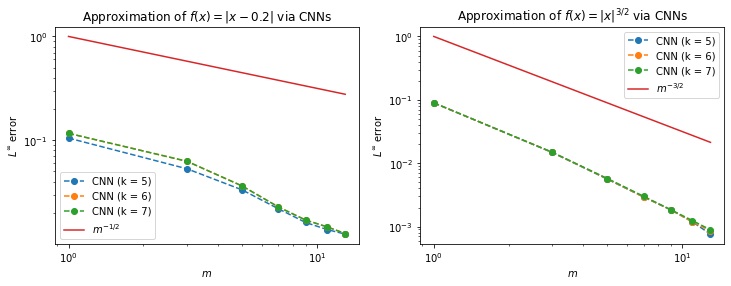}
    \caption{Numerical validation of the upper bounds in Theorem 1. The two panels show the results obtained for signals of different smoothness, respectively $H^{1}(0,1)$ on the left and $H^{2}(0,1)$ on the right. The results are reported for different grid resolutions, the mesh stepsize being $h=2^{-k}$.}
    \label{fig:fdecay}
\end{figure*}
The pictures in Figure \ref{fig:fdecay} show the approximation rates obtained by the actual implementation of $\Psi$ (and $T$) along the lines detailed in the proof. The code used to obtain these results -- as well as all the other results reported in the paper -- is written in Python3 by exploiting the Pytorch library for CNNs, and it is available upon request. Note that we do not train the network, as we directly initialize $\Psi$ with the wished weights and biases. The left panel in Figure \ref{fig:fdecay} shows the results obtained for a mildly smooth signal, $f(x)=|x-1/5|$. In this case we have $s=1$, and the $L^{\infty}$ error between the desired output, $f$, and the CNN approximation, is shown to decay at the expected rate, that is $1/\sqrt{m}$. This is also true regardless of the grid resolution, coherently with Theorem \ref{theorem:single}. Indeed, 
we obtained nearly the same results for $N_{h}=33,65,129$. Finally, the right panel in Figure \ref{fig:fdecay} refers to a smoother case, $f(x)=|x|^{3/2}$, where $s=2$. Here we can remark, once again, the expected behavior.

\section{Learning operators in a parametrized setting}
\label{sec:param}

We now extend the results in Section \ref{sec:approx} to the parameter dependent setting. Let $\Theta\subset\mathbb{R}^{p}$ be some compact parameter space. We are interested in the approximation of an operator $$\Theta\ni\mup\to u_{\mup}\in H^{s}(\Omega).$$ This framework is typically encountered in the case of parameter dependent PDEs, where each value of the input parameter vector $\mup\in\Theta$ usually defines a different PDE solution $u_{\mup}$. In this case, the approximation of the so-called parameter-to-solution map $\mup\to u_{\mup}$ is of remarkable importance, especially when it comes to expensive many-query routines. For the interested reader, we refer to the general literature on Reduced Order Modeling for PDEs \cite{hesthaven,quarteroni} and to more recent contributions on the use of DNNs for the nonintrusive construction of efficient reduced order models in this context \cite{franco2021deep,fresca2021comprehensive,lee2020model}.
\\\\
Here, we aim at characterizing the approximation of such an operator in terms of convolutional neural networks. In particular, we seek for some DNN $\Phi:\mathbb{R}^{p}\to\mathbb{R}^{N_{h}}$ such that
$$\Phi_{i}(\mup)\approx u_{\mup}(x_{i}),$$
where the nodes $\{x_{i}\}_{i}^{N_{h}}\subset\overline{\Omega}$ come from a given discretization of the domain. We build $\Phi$ by considering an architecture that is made of two  blocks, $\phi$ and $\Psi$. The former consists of dense layers, and it has the purpose of pre-processing the input. The latter is instead of convolutional type, and it is used to provide the desired output. We design $\Psi$ along the lines of Theorem \ref{theorem:single}, thanks to which we are able to characterize the approximation error in terms of the network architecture as a whole. 

In particular, we show that: (i) the depth of the dense block depends on the desired accuracy, while the number of convolutional layers only depends on the chosen discretization, (ii) fewer channels in the convolutional layers are required to approximate operators that have highly regular outputs, (iii) the width of the dense layers depends on the smoothness of the operator. We formalize these statements in the Theorem below.
\newpage
\begin{theorem}
\label{theorem:param}
Let $\Omega:=(0,1)$ and let $\{x_{j}\}_{j=1}^{N_{h}}\subset\Omega$ be a uniform grid with stepsize $h=2^{-k}$. Let  $\Theta\subset\mathbb{R}^{p}\ni\mup\to u_{\mup}\in H^{s}(\Omega)$ be  a (nonlinear) operator, where $\Theta$ is a compact domain and $s\ge1$. For some $r\ge0$, assume that the operator is \review{$r$-times continuously Fréchet differentiable. For any $0<\varepsilon<1/2$,} there exists a deep neural network $\Phi:\mathbb{R}^{p}\to\mathbb{R}^{N_{h}}$ such that
$$|u_{\mup}(x_{j})-\Phi_{j}(\mup)|<\varepsilon$$
uniformly for all $\mup\in\Theta$ and all $j=1,\dots,N_{h}$. \review{Additionally, $\Phi$ can be defined as the composition of a fully connected network $\phi$ and a convolutional neural network $\Psi$, i.e. $\Phi=\Psi\circ\phi$, such that the overall architecture has at most}
\begin{itemize}
    \item [i)] $C\log(1/\varepsilon)$ dense layers, with ReLU activation, and $C\log(1/h)$ convolutional layers,
    \item [ii)] $C\varepsilon^{-2/(2s-1)}\left[\review{\varepsilon^{-p/r}}\log(1/\varepsilon)+\log(1/h)\right]$ active\\weights,
    \item [iii)] $C\varepsilon^{-2/(2s-1)}$ channels in input and output,
\end{itemize}
where $C>0$ is some constant dependent on $\Theta$ and on the operator $\mup\to u_{\mup}$, thus also on $s,r,p$.
\end{theorem}
\noindent\textit{Proof}. Let $\varepsilon>0$ and let $c=c(s)>0$ be the constant in Theorem \ref{theorem:single}. We take advantage of the compactness of $\Theta$ and the continuity of the operator to define
$$M:=\max_{\mup\in\Theta}||u_{\mup}||_{H^{s}(\Omega)}<+\infty.$$
Let now $\Psi$ be the CNN in Theorem \ref{theorem:single}, where we fix $m = \lceil (\varepsilon/2)^{-2/(2s-1)}Mc\rceil$. Then, 
\begin{equation}
    \tag{*}
    |u_{\mup}(x_{j})-\Psi_{j}(Tu_{\mup})|< \varepsilon/2
\end{equation}
for all $j=1,\dots,N_{h}$ and $\mup\in\Theta$, where $T: H^{s}(\Omega)\to\mathbb{C}^{2m+1}\cong\mathbb{R}^{4m+2}$ is some continuous linear operator. We now note that, by composition, the map $$\mup \to Tu_{\mup}$$ is an element of the Sobolev space $\review{W^{r,\infty}(\Theta;\mathbb{R}^{4m+1})}$. 

In particular, by Theorem 1 in \cite{yarotsky}, there exists a ReLU DNN $\phi:\mathbb{R}^{p}\to\mathbb{R}^{4m+2}$ with $C\log(1/\varepsilon)$ hidden layers and $Cm\review{\varepsilon^{-p/r}}\log(1/\varepsilon)$ active weights, such that
$$\sup_{\mup\in\Theta}||Tu_{\mup} - \phi(\mup)||_{1}<\varepsilon/2,$$
where $||\cdot||_{1}$ is the $\ell_{1}$ norm over $\mathbb{R}^{4m+2}$, while $C>0$ is a constant that depends on $r,p,\Theta,s$ and the operator $\mup\to u_{\mup}$. The dependence on $s$ comes from the Lipschitz constant of $T$, which may inflate the magnitude of the partial derivatives of $\mup\to Tu_{\mup}$. 

Let now consider the composition $\Phi:=\Psi\circ\phi$. It is easy to see that this DNN architecture satisfies the requirements claimed in the Theorem as soon as we replace the constant $C$ with $\tilde{C}:=Cc$. Also, for any $\mup\in\Theta$ and $j=1,\dots,N_{h}$ we have the desired bound. In fact, by (*),
\begin{multline*}
    \tag{**}
    |u_{\mup}(x_{j})-\Phi_{j}(\mup)|\\\le
    |u_{\mup}(x_{j})-\Psi_{j}(Tu_{\mup})|+|\Phi_{j}(\mup)-\Psi_{j}(Tu_{\mup})|\\<
    \frac{\varepsilon}{2} + |\Psi_{j}(\phi(\mup))-\Psi_{j}(Tu_{\mup})|.
\end{multline*}
Now, we note that $|\Psi_{j}(\mathbf{a})-\Psi_{j}(\mathbf{b})|\le||\mathbf{a}-\mathbf{b}||_{1}$. In fact, in Theorem 1, $\Psi$ was defined as $E\circ R\circ \mathcal{S}_{m}$, where $E$ and $R$ were reshape layers, while $\mathcal{S}_{m}$ was as in Lemma \ref{lemma:fourier}. In particular,
$$|\Psi_{j}(\mathbf{a})-\Psi_{j}(\mathbf{b})|\le \sup_{x\in[0,1]}\left|\sum_{k=-m}^{m}a_{k}\textrm{e}^{2\pi\I kx}-\sum_{k=-m}^{m}b_{k}\textrm{e}^{2\pi\I kx}\right|$$$$\hspace{16em}\le\sum_{k=-m}^{m}|a_{k}-b_{k}|.$$
Therefore, relationship (**) finally yields
$$|u_{\mup}(x_{j})-\Phi_{j}(\mup)|\le\dots<\frac{\varepsilon}{2} + ||\phi(\mup)-Tu_{\mup}||_{1}<\varepsilon.$$
$\proved$

\theoremstyle{remark}
\newtheorem{remark}{Remark}
\begin{remark}
\label{remark}
\normalfont
\review{The hypothesis of Frechét differentiability in Theorem \ref{theorem:param} can be relaxed. Indeed, for the error bounds to hold, it is sufficient that the map $\mup\to Tu_{\mup}$ is in the Sobolev space $W^{r,+\infty}(\Theta,\mathbb{R}^{4m+1})$. For instance, one may require the operator $\mup\to u_{\mup}$ to be $r$-times Frechét differentiable, with the first $r-1$ derivatives being continuous and the last one being (essentially) bounded.
}
\end{remark}

\begin{remark}
\label{remark2}
\normalfont
\review{Ultimately, the proof of Theorem \ref{theorem:param} is based on a representation of the form $u_{\mup}(\x)\approx\sum a_{i}(\mup)v_{i}(\x)$, where the sum runs over $\mathcal{O}(m)$ terms. The dense block learns the map $T:\mup \to [a_{i}(\mup)]_{i}$, while the convolutional block expands the coefficients over the spatial grid, implicitly reconstructing the $v_{i}$ modes. However, we remark that the latter are not precisely the Fourier modes and, similarly, the $a_{i}$ are not the Fourier coefficients of $u_{\mup}$. For instance, as detailed in the proof of Theorem \ref{theorem:single}, the parameter-to-coefficients operator, $T$, includes some additional artifacts, such as an implicit (smooth) periodicization of the signals. All these arguments are fundamental for the proof of Theorem \ref{theorem:param} to work out, but they are not required in practical applications. We do not need to specify or know the operator $T$: the dense block will automatically find a suitable representation from data, while the convolutional part will learn a corresponding basis expansion.
In this sense, there is a clear analogy between our construction and the tools that other researchers have been using for the study of DeepONets: see, e.g., the work by \cite{de2022generic}.
}
\end{remark}

\begin{remark}
    \lastchanges{It should be noted that, with respect to the size of the input dimension, $p$, the result in Theorem \ref{theorem:param} suffers from the curse of dimensionality. In general, however, one can only achieve better approximation rates by either including additional assumptions about the operator (for instance, by resorting to the case of $p$-holomorphic functions, as in \cite{schwab2019deep}) or by considering different metrics when evaluating errors (e.g., by moving from a worst-case approach to a probabilistic one, as in \cite{lanthaler2022error}). As for now, the usefulness of our results is mostly limited to the situation in which $p\ll N_{h}$, which, in any case, is a scenario fairly common to many applications.}
\end{remark}

\section{Numerical validation}
\label{sec:valid}
We finally present some numerical experiments that confirm the decay rates predicted in Theorem \ref{theorem:param}. We proceed as follows. Once introduced the operator to be learned, we identify the smoothness indices $s$ and $r$ that appear in Theorem \ref{theorem:param}. Then, we fix a guess architecture $\Phi^{(1)}$ that serves as a starting point. Following the ideas of Theorem \ref{theorem:param}, we prescribe $\Phi^{(1)}$ as a DNN that is made by two blocks, $$\Phi^{(1)}=\Psi^{(1)}\circ\phi^{(1)},$$ where $\phi^{(1)}$ is dense, while $\Psi^{(1)}$ is of convolutional type. More precisely:
\review{
\begin{itemize}
    \item we let $\phi^{(1)}$ have $L_{1}$ hidden layers of constant width $w_{1}$, while we equip the output layer with $m_{1}(2^\ell+5)$ neurons. In this way, the output of $\phi^{(i)}$ can be reshaped in the form of a $m_{1}\times(2^\ell+5)$ matrix, which allows the CNN block to interpret it as a discrete signal having $m_{1}$ features of length $2^\ell+5$;
    \item we let $\Psi^{(1)}:=R_{2}\circ\tilde{\Psi}^{(1)}\circ R_{1}$, where the $R_{i}$ are auxiliary reshape operations, whereas $\tilde{\Psi}^{(1)}$ contains the actual convolutional part. In particular, $R_{1}$ is used to reshape the output of $\phi^{(1)}$ as $m_{1}\times(2^\ell+5)$. Conversely, $\tilde{\Psi}^{(1)}$ is comprised of $k-\ell+1$ transposed convolutional layers, where $k:=\log_{2}(1/h)$ depends on the final grid resolution. All these layers have $m_{1}$ channels at input and output, grouped by $m_{1}$ (cf. Definition \ref{def:cnnt} in the Appendix). The architecture of $\tilde{\Psi}^{(1)}$ is then further supplemented with a convolutional layer having a single output channel, and a terminal transposed convolution. All layers in $\tilde{\Psi}^{(1)}$ use a kernel of size $5$ and a stride of $2$, except for the last one, which resorts to the default stride of $1$. 
    With this set up, the output of $\tilde{\Psi}^{(1)}\circ R_{1}\circ \phi^{(1)}$ is guaranteed to have shape $1\times (2^{k+1}+1)$. The final reshaping, $R_{2}$, is then used to flatten and half the output, leaving us with the correct dimension, that is, $2^{k}+1$.
\end{itemize}}
\newcommand{\ntrain}{N_{\text{train}}}\newcommand{\fomdim}{N_{h}}\newcommand{\test}{\text{test}}\newcommand{\ntest}{N_\test}
\review{Table \ref{tab:architectures} reports in full detail the complete structure of the neural network architecture. We remark that the latter is itself parametrized by five hyperparameters, namely $L=L_{1},w=w_{1},m=m_{1},k$ and $\ell$. The first three allow us to tune the overall complexity of the model, and we shall exploit them to verify the estimates in Theorem \ref{theorem:param}. Conversely, the last two, $k$ and $\ell$, are problem dependent and we shall fix their value for each experiment alone. In fact, $k$ is related to the grid discretization, while $\ell$ is used to define an intermediate level of resolution.}
\\
\\
We initialize all the weights and biases of $\Phi^{(1)}$ randomly, following the approach introduced by He et al. in \cite{he}. We then train $\Phi^{(1)}$ over a training set $\{\mup_{i}, u_{\mup_i}\}_{i=1}^{\ntrain}$ in such a way that the loss function below is minimized
\begin{equation}
    \label{eq:loss}
    \mathcal{L}(\Phi^{(1)}):=\frac{1}{\ntrain}\sum_{i=1}^{\ntrain}\left(h\sum_{j=1}^{\fomdim}|u_{\mup_i}(x_{j})-\Phi^{(1)}_{j}(\mup_{i})|^{2}\right),
\end{equation}where $x_{1},\dots,x_{N_{h}}$ is some dyadic partition of $(0,1)$ associated to a given grid resolution $h=2^{-k}$.
We then evaluate $\Phi^{(1)}$ over a test set of unseen instances $\{\mup^\test_{i}, u_{\mup^\test_i}\}_{i=1}^{\ntest}$ in order to compute the empirical \review{uniform} error, given by
\begin{equation}
\label{eq:test}
E(\Phi^{(1)}) = \max_{i,j}|u_{\mup^{\test}_i}(x_{j})-\Phi^{(1)}_{j}(\mup^{\test}_{i})|.    
\end{equation}
Then, we exploit Theorem \ref{theorem:param} in an attempt to define a second architecture, $\Phi^{(2)}$, that can be twice as accurate by halving the error over the testing set, that is by requiring that $E(\Phi^{(2)})\approx E(\Phi^{(1)})/2$.
We do this as follows:
\begin{itemize}
    \item we update the number of channels according to (iii) in Theorem \ref{theorem:param}. In particular, up to rounding operations, we let 
    $$m_{2}:=2^{2/(2s-1)}m_{1};$$
    \item we increase the with of the dense layers coherently with (ii) in Theorem \ref{theorem:param}, that is
    $$w_{2}:=\review{\sqrt{2^{p/r+2/(2s-1)}}}w_{1},$$
    \review{where the square root comes from the fact that a dense layer from $\mathbb{R}^{w}\to\mathbb{R}^{w}$ carries $\mathcal{O}(w^{2})$ active weights};
     \item as suggested by (i) in Theorem \ref{theorem:param}, we also increase the number of dense hidden layers. In principle, the depth of the dense block should be increased by a constant factor $C\log(2)$. In practice, we let
     $$L_{2} = L_{1} + l,$$
     where $l$ is either 1 or 2. This is to ensure that the obtained architectures are still feasible to train, as very deep models may become hard and expensive to optimize.
\end{itemize}
We then train $\Phi^{(2)}$ and iterate the above steps to generate $\Phi^{(3)}$, so that $$E(\Phi^{(j)})\propto 2^{-j}.$$
We highlight that, according to Theorem \ref{theorem:param}, this procedure should be robust with respect to the space discretization. In other words, we expect to obtain similar results regardless of the number of grid points employed in the discretization. To assess whether this behavior is actually observed in practice, we repeat our analysis for different mesh step sizes $h=2^{-k}$ (when possible).
 
\renewcommand{\arraystretch}{1.5}
\begin{table}
\centering
    \begin{tabular}{|p{1.5cm}|p{1.75cm}|p{1.75cm}|p{1.75cm}|}
    \hline
    \textbf{Layer} & \textbf{Input} & \textbf{Output} & \textbf{Active weights}
    \\
    \hline\hline
    Dense & $p$ & $w$  & $pw$ \\
    Dense & $w$ & $w$  & $w^{2}$ \\
    \end{tabular}
    \begin{center}...\\\textit{Repeat last layer L-1 times}\\...\end{center}
    \begin{tabular}{|p{1.5cm}|p{1.75cm}|p{1.75cm}|p{1.75cm}|}
    Dense &  $w$ & $(2^\ell +5)m$  & $(2^\ell+5) wm$ \\
    Reshape &  $(2^\ell+5) m$ & $m\times(2^\ell+5)$  & - \\
    ConvTr$_{5,2}$ & $m\times(2^\ell+5)$ & $m\times(2^{\ell+1}+13)$  & 5$m$ \\
  \end{tabular}
    \begin{center}...\\\textit{Repeat last layer $k-\ell+2$ times}\\...\end{center}
    \begin{tabular}{|p{1.5cm}|p{1.75cm}|p{1.75cm}|p{1.75cm}|}
    Conv$_{5,2}$  & $m\times(2^{k+2}-3)$ & $1\times(2^{k+1}-3)$  & 5$m$ \\
    ConvTr$_{5,1}$  & $1\times(2^{k+1}-3)$ & $1\times(2^{k+1}+1)$  & 5 \\
    Reshape  & $1\times(2^{k+1}+1)$ & $2^{k+1}+1$  & 0 \\
    Truncate &  $2^{k+1}+1$ & $2^{k}+1$  & 0\\
    \hline
    \end{tabular}\vspace{0.5em}
    \caption{\review{Parametric architecture for the numerical experiments, Section \ref{sec:valid}. The dense block depends on the structural hyperparameters $p$ (input dimension), $w$ (width of the hidden layers), $L$ (number of hidden layers-1), $m$ and $\ell$ (intermediate output). Conversely, the design of the convolutional part depends on $m$ (number of channels), $\ell$ and $k=\log_{2}(1/h)$ (depth of the convolutional block, upto a constant), where $h$ is the stepsize of the discretization. Here, Conv$_{s,t}$ = Convolutional layer with kernel of size $s$ and stride $t$; ConvTr$_{s,t}$ = Transposed convolutional layer with kernel of size $s$ and stride $t$; Truncate = weightless layer that keeps only the central $2^{k}+1$ components of the input. All learnable layers, except the last one, employ the 0.1-leakyReLU activation.}}
    \label{tab:architectures}
\end{table}
 
\subsection{Benchmark example}
\label{bench}
To start, we consider the approximation of an operator that is defined analytically. More precisely, let $\Theta = [0,1]\times[0,1]\times[1,2]\subset\mathbb{R}^{3}$. For any fixed $\mup=(\mu_{1}, \mu_{2}, \mu_{3})\in\Theta$ let
$$u_{\mup}(x) = \mu_{3}|x-\mu_{1}|^{3}\textrm{e}^{-\mu_{2}x}.$$
We are interested in learning the map $\mup\to u_{\mup}$. To this end, we note that $\{u_{\mup}\}_{\mup\in\Theta}\subset H^{3}(\Omega)\setminus H^{4}(\Omega)$. Also, the operator is at most twice differentiable with respect to $\mup$, as its third derivative becomes discontinuous. According to the notation in Theorem \ref{theorem:param}, this results in $s=3$ and $r=2$. \review{However, due to the boundness of the third derivative, we can actually apply Theorem \ref{theorem:param} with an increased smoothness index, i.e. $r=3$ (see Remark \ref{remark}).}

\begin{figure}[tb]
    \centering
    \includegraphics[width=0.65\textwidth]{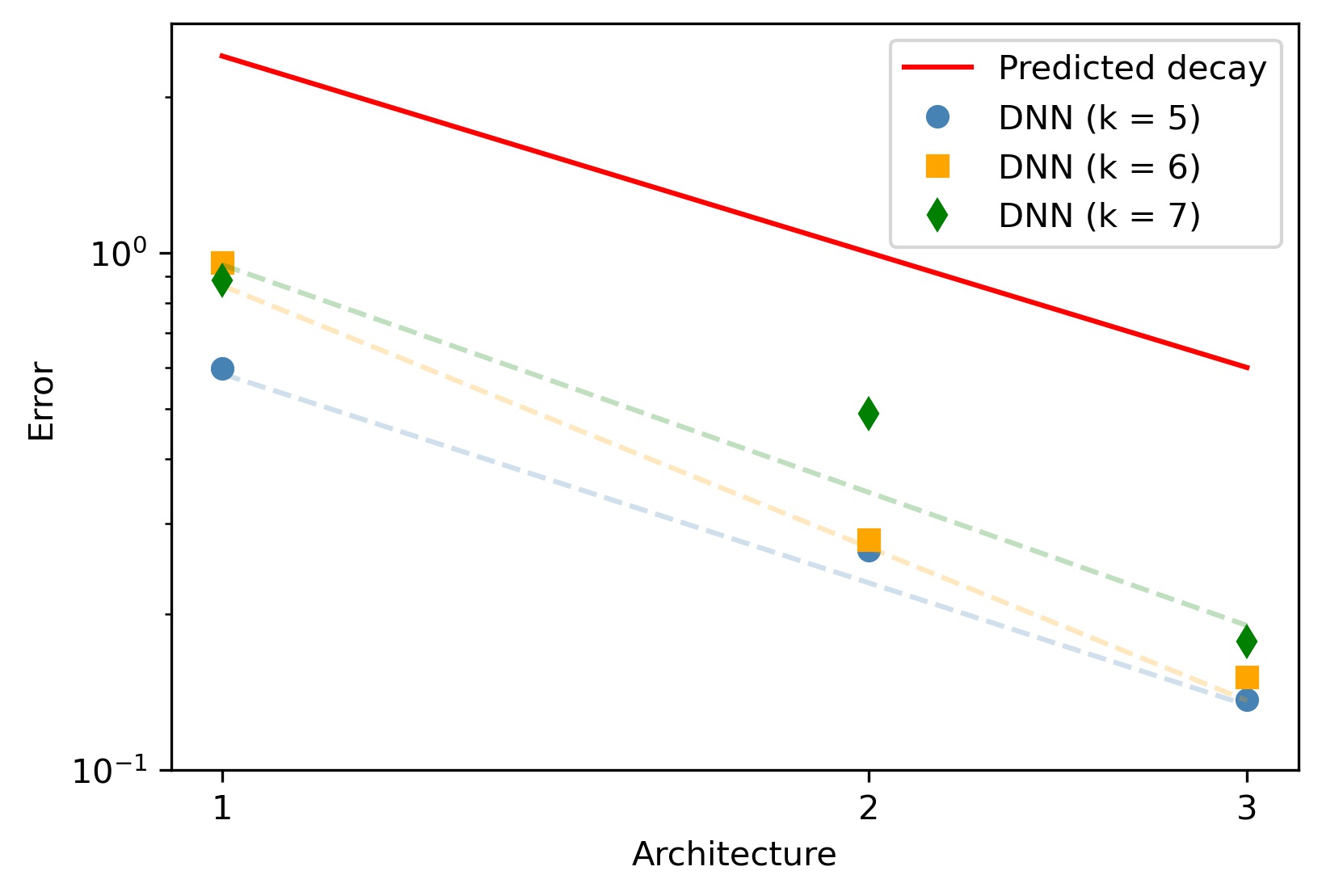}
    \caption{Numerical validation of Theorem \ref{theorem:param} for the Benchmark example, Section \ref{bench}. Axis are in loglog scale. The red line corresponds to the predicted decay rate, $2^{-j}$, while the markers refer to the DNN models. Different markers correspond to different grid resolutions. The errors $E(\Phi^{(j)})$ obtained for the architectures $j=1,2,3,$ are computed as in Equation \eqref{eq:test}. \review{Dashed lines are obtain through least square estimates.}}
    \label{fig:paramdecay}
\end{figure}

For the space discretization, we consider three different mesh resolutions, $h=2^{-5},2^{-6},2^{-7}$, corresponding respectively to $N_{h}=33,65,129$ grid points. We train the networks by minimizing \eqref{eq:loss} via the so-called L-BFGS optimizer, where the training set consists of 500 randomly sampled parameters instances. We do not use batching strategies and we set the learning rate to its default value of 1. To avoid possible biases introduced by the optimization, we initialize and train each architecture multiple times (here, five), only to keep the best out of all the training sessions. This is a common practice known as \textit{ensemble training}. For our starting architecture, $\Phi^{(1)}$, we set \review{$$m_{1}=5,\;w_{1}=1,\;L_{1}=1.$$} In this case, we have $2^{2/(2s-1)}\approx1.32$ and $\review{\sqrt{2^{p/r+2/(2s-1)}}}\approx1.62,$ since $p=3$. In particular, our strategy for enriching the architectures can be stated as follows: to obtain a model that is twice as accurate, we increase the number of channels in the convolutional layers by nearly 30\%, \review{while we add about 50\% new neurons to the dense layers}. The theory also suggests to increase the depth of the dense block by some constant factor $l$. Here, we let $l=1$. \review{Finally, for this numerical experiment we choose a fixed coarse resolution of $13=2^{\ell}+5$, that is, we let $\ell=3$ in Table \ref{tab:architectures}.}
\begin{figure*}
    \centering
    \includegraphics[width=\textwidth]{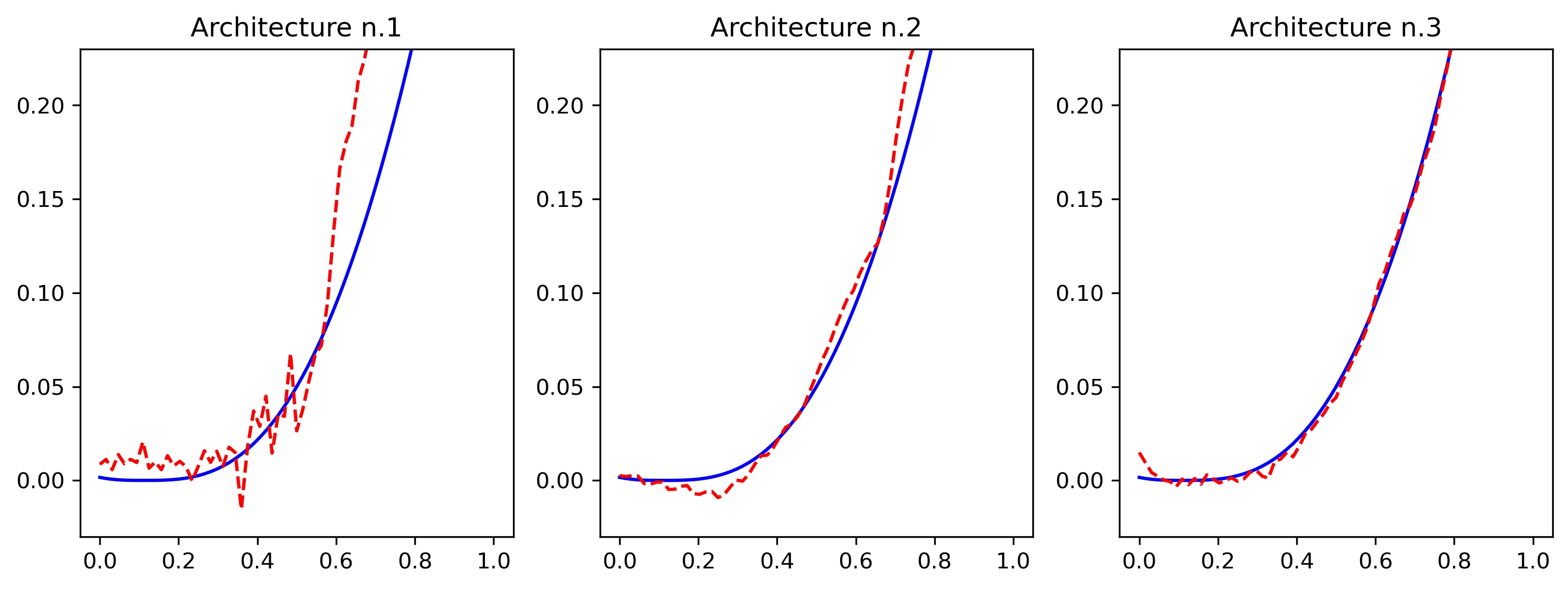}
    \caption{Benchmark example, Section \ref{bench}. \review{In blue, an instance $u_{\mup}$ coming from the test set, here for $\mup=[0.1125, 0.5409, 0.1255].$ The red dashed lines report the approximations proposed by the three DNNs, respectively $\Phi^{(1)},\Phi^{(2)},\Phi^{(3)}$. Grid resolution is $h=2^{-6}$.}}
    \label{fig:paramex}
\end{figure*}\\\\
Results are in Table \ref{tab:param1}, Figures \ref{fig:paramex} and \ref{fig:paramdecay}. The first picture compares the output of the three architectures with that of the operator, for an unseen value of the input parameter $\mup$. The quality of the approximation clearly increases as we consider richer and richer models. \review{In general, we see that the estimated signals are rougher compared to the ground truth. This, however, is most likely due to our use of the leaky ReLU activation (which we chose in order to be consistent with Theorem \ref{theorem:param}): other activations may lead to smoother results, with a possible benefit in terms of approximation properties, see e.g. \cite{guhring2021approximation}, or training strategies, see e.g. \cite{mishra2021enhancing}.}
We also note that the regions with lower regularity are the most difficult to capture, coherently with what we expected.\renewcommand{\arraystretch}{1.5}
\begin{table}
    \centering
    \begin{tabular}{|c|c|c|c|c|l|c|}
    \hline
    Model & $m_{j}$   &  $w_{j}$ & $L_{j}$ & Active weights & $E(\Phi^{(j)})$ \\
    \hline\hline
    $\Phi^{(1)}$ &  5 & 1 & 1  & 223 & 0.884 \\
    $\Phi^{(2)}$ &  7 & 2 & 2 & 407 & 0.489 \\
    $\Phi^{(3)}$ &  9&  3 & 3 & 653 & 0.177 \\
       \hline
    \end{tabular}\vspace{0.5em}
    \caption{\review{Architectures and corresponding errors for the Benchmark example, Section \ref{bench}. Results are reported limitedly to the case of grid resolution $h=2^{-7}$. The hyperparameters read as in Section \ref{sec:valid}, that is: $m_{j}$ = maximum number of convolutional features in the CNN block, upto a multiplicative constant; $w_{j}$ = number of neurons per dense layer; $L_{j}$  = depth of the dense block. The errors $E(\Phi^{(j)})$ are computed as in Equation \eqref{eq:test}}.}
    \label{tab:param1}
\end{table}
 \review{Indeed,
the architectures mostly struggle in capturing flat regions, which is understandable as these entail discontinuities in the higher-derivatives.
Finally,} Figure \ref{fig:paramdecay}, reports the errors $E(\Phi^{(j)})$ in comparison with the expected decay rate $2^{-j}$. \review{There,} we see that the numerical results perfectly match the theory, regardless of the grid resolution that is considered.

\subsection{Application to a parametrized time-dependent nonlinear PDE}
\label{pde}
We now consider a benchmark consisting of a one-dimensional coupled PDE-ODE nonlinear system
\begin{equation}
\label{eq:412-1}
\begin{cases}
\vspace{0.1cm}
\displaystyle \mu \frac{\partial u_{\mup}}{\partial t} - \mu^2 \frac{\partial^2 u_{\mup}}{\partial x^2} + R(u_{\mup}) + w_{\mup} = 0, \quad & (x, t) \in \Omega \times (0,T) \\
\vspace{0.1cm}
\displaystyle \frac{d w_{\mup}}{d t} + (2 w_{\mup} - 0.5 u_{\mup})=0, \quad & (x, t) \in \Omega \times (0,T) \\
\vspace{0.1cm}
\displaystyle \frac{\partial u_{\mup}}{\partial x}(0,t) = 50000 t^3 e^{-15t}, \quad & t \in (0,T)\\
\vspace{0.1cm}
\displaystyle \frac{\partial u_{\mup}}{\partial x}(1,t) = 0, \quad & t \in (0,T)\\
u_{\mup}(x,0)=0, \; w_{\mup}(x,0)= 0, \quad & x \in \Omega,
\end{cases}
\end{equation}
where $R(u_{\mup}):=u_{\mup}(u_{\mup}-0.1)(u_{\mup}-1)$, while $\Omega = (0, 1)$ and $T = 2$. The above consists in a parametrized version of the monodomain equation coupled with the FitzHugh-Nagumo cellular model, describing  the excitation-relaxation of the cell membrane in the cardiac tissue \cite{fitzhugh1961impulses,nagumo1962anactive}. System \eqref{eq:412-1} has been discretized in space through linear finite elements, by considering $N_h = 2^k$, with $k \in \mathbb{N}$, grid points, and using a one-step, semi-implicit, first order scheme  for time discretization with time-step $\Delta t = 5 \times 10^{-3}$;  see, e.g., \cite{pagani2018numerical} for further details.  
The solution of the former problem consists in a parameter-dependent traveling wave, which exhibits sharper and sharper fronts as the parameter $\mu$ gets smaller. The numerical transmembrane potential solution $u_{\mup}$ represent the ground truth data in the experiments reported in the following.
\\\\
Here, we consider the map $(\mu, t)\to u_{\mup}(\cdot, t)$ as our operator of interest. In particular, the two dimensional vector parameter $\mup:=(\mu,t)$ consists of the scalar parameter $\mu$ and the time variable $t$. We let $\mup$ vary in the (time-extended) parameter space $\Theta:=\Theta_{0}\times[0, T]$, where $\Theta_{0} := 5 \cdot  [10^{-3}, 10^{-2}]$.

In this case, it is not straightforward to identify the smoothness indices $s$ and $r$. The numerical simulations show that the solutions $u_{\mup}$ to \eqref{eq:412-1} tend to have sharp gradients for certain values of the scalar parameter $\mu$. In light of this, we let $s=1$; if the solutions are actually smoother, then we expect the errors to decay faster than the predicted rate. Conversely, we make the assumption that $r = +\infty$, i.e. that the parameter-to-solution map is infinitely differentiable. We remark that the constant $C$ appearing in Theorem \ref{theorem:param} actually depends on $r$. To this end, we make the further assumption that $C=C(r)$ is bounded with respect to $r$. \review{This is a rather restrictive assumption, but we expect the latter to hold for analytic operators with fast decaying coefficients: indeed, in this case, one could adapt the proof of Theorem \ref{theorem:param} by replacing the result due to Yarotsky with a stronger one, such as, e.g., Theorem 3.9 in  \cite{schwab2019deep}}.

As a starting point, we consider the following structural hyperparameters 
\review{$$m_{1}=1,\;w_{1}=5,\;L_{1}=1,$$}
to build our reference architecture $\Phi^{(1)}$. In this case, we have $2^{\review{p/r}+2/(2s-1)}=2^{2/(2s-1)}=4.$ In particular, in order to half the test error, Theorem \ref{theorem:param} suggests to quadruplicate the number channels in the convolutional layers, and to double that of the neurons in the dense layers. Regarding the depth of the dense block, instead, we increase it by a constant factor of $l=1$ when moving from an architecture to a more complex one. 
In this case, we do not assess the model performance for varying resolution levels as we stick to the same grid employed by the Finite Element solver. Instead, in order to collect more data, we repeat the same analysis for a different guess architecture, namely
\review{$$m_{1}=1,\;w_{1}=2,\;L_{1}=2.$$}

\begin{figure*}
    \centering
    \includegraphics[width=\textwidth]{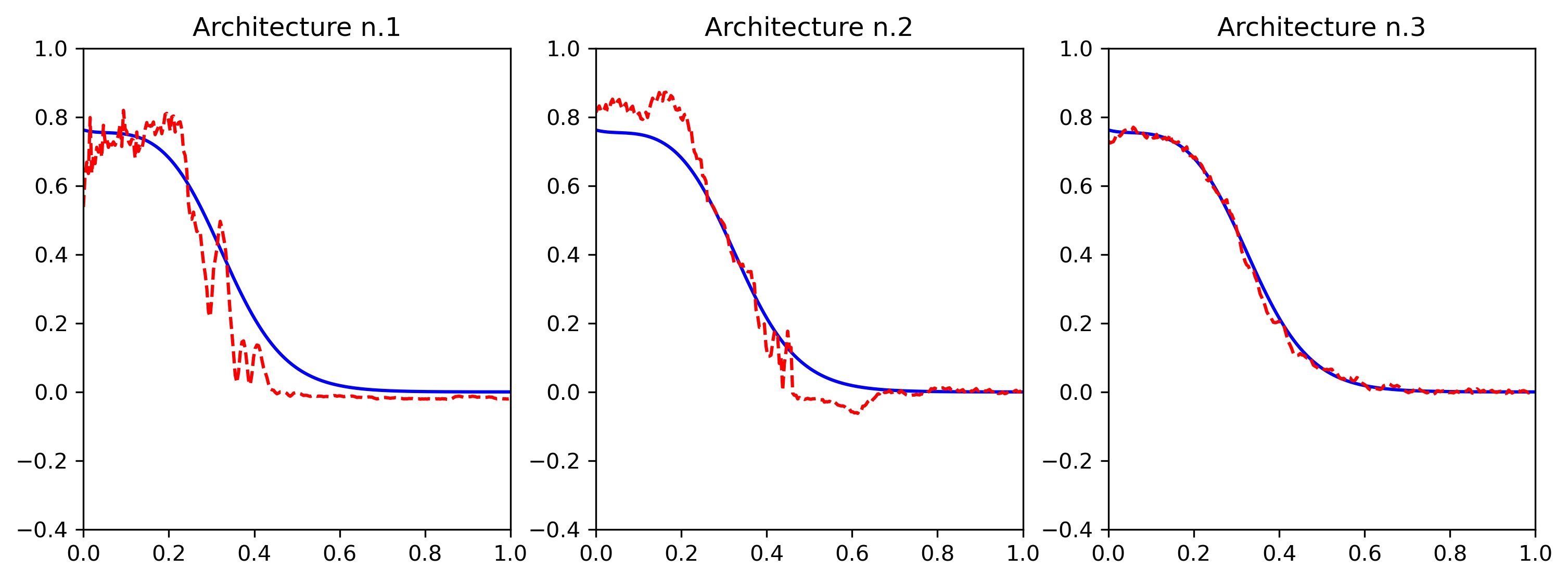}
    \caption{Learning the parameter-to-solution operator of a parametrized time-dependent nonlinear PDE, Section \ref{pde}. In blue, a snapshot $u_{\mup}(\cdot, t)$ coming from the test set, here for $\mu=0.0488$ and  $t= 0.7250.$ The red dashed lines correspond to the approximations proposed by the three DNN models, respectively $\Phi^{(1)},\Phi^{(2)},\Phi^{(3)}$.}
    \label{fig:monoex}
\end{figure*}

To collect the training and test sets, we proceed as follows. We sample $N_{train} = 20$ equally spaced values for the scalar parameter $\mu\in\Theta_{0}$, and we consider their midpoints to obtain $N_{test} = 19$ test instances. For each $\mu\in\Theta_{0}$ fixed, we then extract uniformly $N_t = 25$ time snapshots from the global trajectory defined over the interval $[0,T]$. Once again, we train the DNN models using the L-BFGS optimizer (no batching, learning rate = 1).
\begin{figure}
    \centering
    \includegraphics[width=0.65\textwidth]{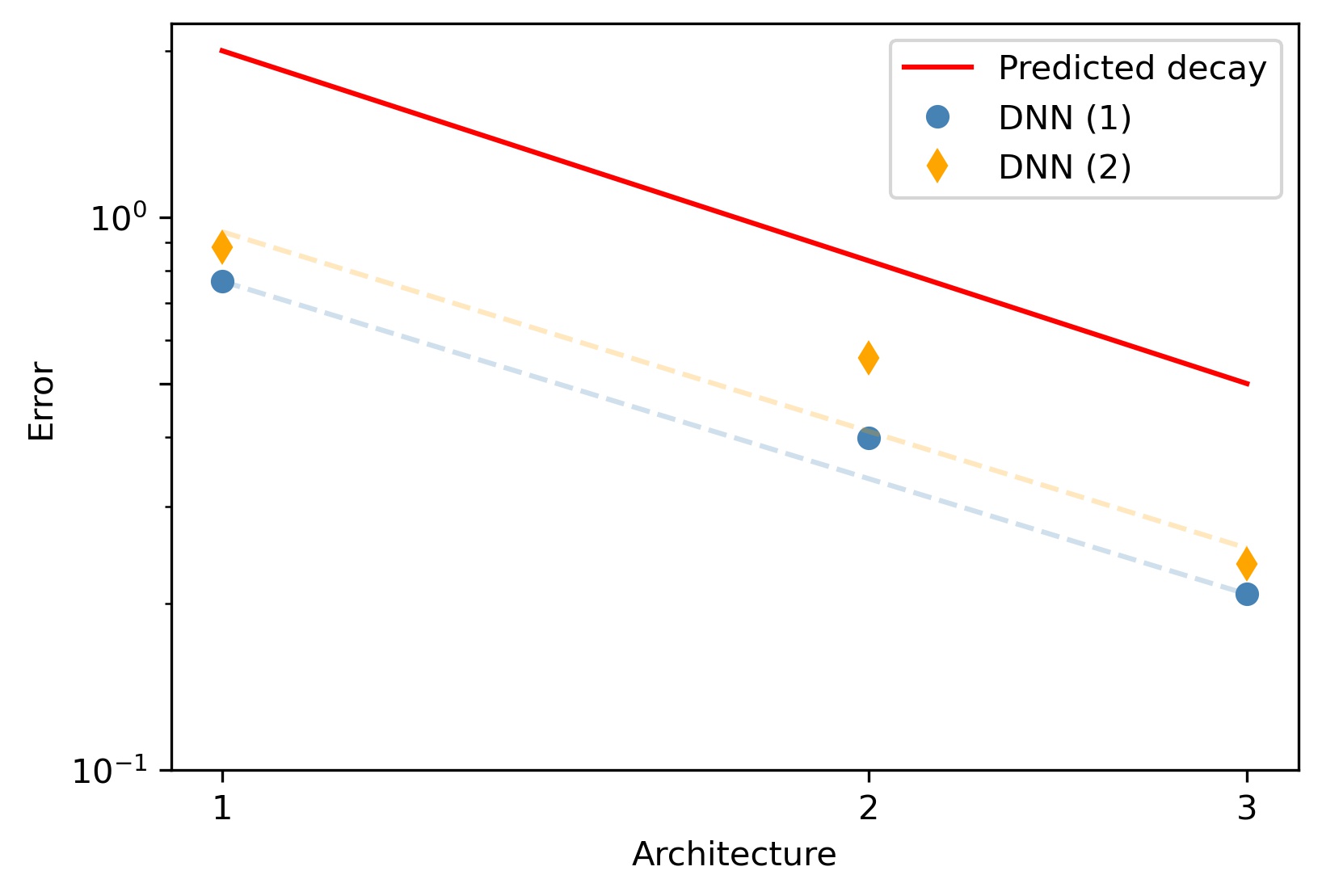}
    \caption{Numerical validation of Theorem \ref{theorem:param} for the coupled FitzHugh-Nagumo problem, Section \ref{pde}. Axis read in loglog scale. The red line corresponds to the predicted decay rate, $2^{-J}$, while the blue markers report the errors $E(\Phi^{(j)})$ obtained for the architectures $j=1,2,3.$ Different markers correspond to different choices of the initial guess architecture $\Phi^{(1)}$, respectively. Errors are computed accordingly to Equation \ref{eq:test}.}
    \label{fig:monodecay}
\end{figure}
\renewcommand{\arraystretch}{1.5}
\begin{table}
    \centering
    \begin{tabular}{|c|c|c|c|l|c|}
    \hline
    Model & $m_{j}$   &  $w_{j}$ & $L_{j}$ & Active weights & $E(\Phi^{(j)})$ \\
    \hline\hline
    $\Phi^{(1)}$ &  1 & 5 & 1  & 225 & 0.765 \\
    $\Phi^{(2)}$ &  4 & 10 & 2 & 1'705 & 0.398 \\
    $\Phi^{(3)}$ &  16 & 20 & 3 & 13'085 & 0.208 \\
       \hline
    \end{tabular}\vspace{0.5em}
    \caption{\review{Architectures and corresponding errors for the coupled FitzHugh-Nagumo problem, Section \ref{pde}. Results are reported limitedly to one of the initial guess architectures. Hyperparameters read as in Section \ref{sec:valid} and Table \ref{tab:param1}. Errors are computed as in Equation \eqref{eq:test}.}}
    \label{tab:param2}
\end{table}
\\\\
Results are reported in Figures \ref{fig:monodecay} and \ref{fig:monoex}. As for the benchmark example, we see that the DNN models become more and more expressive as we move from $\Phi^{(1)}$ to $\Phi^{(3)}$. \review{Furthermore,} the error trend, reported in Figure \ref{fig:monodecay}, is in agreement with the estimates presented in Theorem \ref{theorem:param} regardless of the initial guess for the architecture. Note, once again, that here we only consider one resolution level, as we employ the same step size $h$ adopted by the Finite Element solver.

\begin{figure*}[tb]
    \centering
    \includegraphics[width=\textwidth]{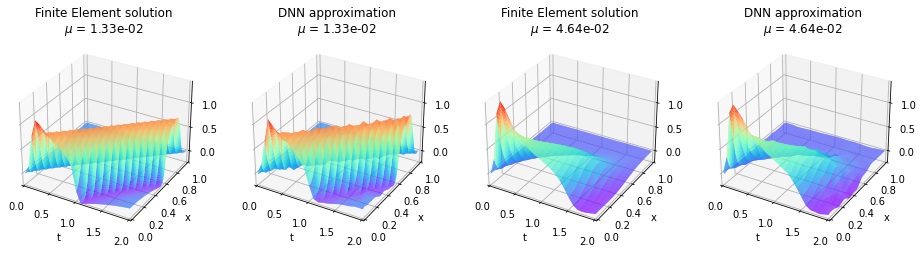}
    \caption{Learning the parameter-to-solution operator of a parametrized time-dependent nonlinear PDE, Section \ref{pde}. Comparison between Finite Element solutions and DNN approximations for different $\mu$. The first and the third plot report the spacetime surface $[0,T]\times\Omega\to\mathbb{R}$ representing the Finite Element simulation, thus $(t,x)\to u_{\mu}(x,t)$. Conversely, the second and the fourth picture show the corresponding DNN approximation over the same spatial grid, $(t,x_{j})\to \Phi_{j}(t,\mu)$. Here, $\Phi$ is constructed considering the third architecture generated during the iterative augmentation process in Section \ref{sec:valid}, starting from the second guess architecture in Section \ref{pde}.}
    \label{fig:monotime}
\end{figure*}
Since we included time as an additional parameter, the plots in Figure \ref{fig:monoex} fix both the scalar parameter $\mu$ and the time instant $t$. However, we recall that according to Equation \eqref{eq:test} the model was evaluated in terms of worst-case errors. In particular, the quality of the approximation is guaranteed over the whole time interval $[0,T]$ and for any choice of the scalar parameter $\mu\in\Theta_{0}$. Figure \ref{fig:monotime}, shows the overall dynamics of the solution for two different choices of $\mu$, with a comparison between Finite Element solutions and DNN approximations. Despite containing a few numerical artifacts, we see that the DNN model fully captures the general behavior of the solutions, both in the hyperbolic and diffusive case ($\mu=1.33\cdot10^{-2}$ and $\mu=4.64\cdot10^{-2}$ respectively). Of note, the spurious oscillations in the DNN approximation are in perfect agreement with the errors reported in Table \ref{tab:param2}. Accordingly to Theorem \ref{theorem:param}, these can be removed by considering larger architectures and, possibly, more training data.

\section{Conclusions}
\label{sec:conc}
In this paper, we have established and verified theoretical error bounds for the approximation of nonlinear operators by means of CNNs. Our results shed a light on the role played by convolutional layers and their hyperparameters, such as input-output channels, depth and others. In particular, they show how operator learning problems can be decoupled in two parts: on the one hand, the difficulty in characterizing the dependence with respect to the input parameters; on the other hand, the issue in having to reconstruct complex space-dependent outputs. The presented research is original and timely. 
Indeed, at the best of our knowledge, all the available results on DNNs and operator learning 
do not address the peculiar properties of CNNs, instead they consider classic fully connected architectures. Conversely, those works that focus on CNN models are typically not framed in the context of operator learning. 

Our analysis is limited to the 1-dimensional case, $d=1$, that is when the output of the operator are functions defined over an interval. However, we note that the main ideas underlying our proofs can be extended to higher dimensions with little effort. The critical points are Lemmas \ref{lemma:fundamental}, \ref{lemma:fourier} and Theorem \ref{theorem:single}. For the first two results, one needs to define suitable convolutional layers that are able to advance along different dimensions separately, which can be carried out via 2D and 3D convolutions whenever $d=2,3$. Conversely, Theorem \ref{theorem:single} has to be adapted in a proper way, since it becomes trickier to turn generic maps $f:[0,1]^{d}\to\mathbb{R}$ onto periodic functions. Furthermore, as the spatial dimension $d$ plays an important role in Sobolev inequalities, it may be convenient to replace the original output space $H^{s}(\Omega)$ with other functional spaces, such as $W^{s,\infty}(\Omega)$ or $\mathcal{C}^{s}(\overline{\Omega})$, when addressing the case $d>1$.

Nevertheless, we believe that our results motivate the recent success of CNNs, especially in areas such as Reduced Order Modeling of PDEs. This is because, as shown in Theorem \ref{theorem:param}, smooth outputs are those that are better approximated by CNNs. Solutions to partial differential equations often enjoy regularity properties that make them an appealing area of application for the proposed analysis. This further promotes the practical use of CNNs as well as their theoretical study from a purely mathematical point of view.

\section*{Fundings and acknowledgements}
NRF, PZ and AM have been partially supported by the ERA-NET ERA PerMed / FRRB grant agreement No ERAPERMED2018-244, RADprecise - Personalized radiotherapy: incorporating cellular response to irradiation in personalized treatment planning to minimize radiation toxicity. SF and AM have been partially supported by Fondazione Cariplo, Italy, grant n. 2019-4608 and by the Italian National Group of Scientific Computing (INDAM-GNCS).
The authors also thank Professor Fabio Nobile (EPFL, Lausanne) for the helpful discussion about this work.

\appendix
\section*{Appendix}
We report, in mathematical terms, the formal definition of convolutional layers and CNNs. These correspond to the ones adopted in the literature and reported within the Pytorch documentation. For tensor objects we use the following notation. 
Given $\mathbf{A}\in\mathbb{R}^{n_{1}\times\dots\times n_{d}}$, we write $\mathbf{A}_{i_{1},\dots,i_{p}}$ for the subtensor in $\mathbb{R}^{n_{p+1}\times\dots\times n_{d}}$ obtained by fixing the first $p$ dimensions along the specified axis, where $1\le i_{j}\le n_{j}$. We also adopt the usual abuse of notation for which scalar-valued activation functions operate componentwise on vectors, that is
$$\rho([x_{1},\dots,x_{l}]):=[\rho(x_{1}),\dots,\rho(x_{l})]$$
whenever $\rho:\mathbb{R}\to\mathbb{R}$.

\begin{definition}
Let $m,m's,t,d$ be positive integers and let $g$ be a common divisor of $m$ and $m'$. A 1D Convolutional layer with $m$ input channels, $m'$ output channels, grouping number $g$, kernel size $s$, stride $t$, dilation factor $d$ and activation function $\rho:\mathbb{R}\to\mathbb{R}$, is a map of the form
$$\Phi:\mathbb{R}^{m\times n}\to\mathbb{R}^{m'\times \left\lfloor\frac{n-d(s-1)-1}{t} + 1\right\rfloor}$$
whose action on a given input $\mathbf{X}\in\mathbb{R}^{m\times n}$ is defined as
 $$\Phi(\mathbf{X})_{k'}=\rho\left(\sum_{k}\mathbf{W}_{k',k}\otimes_{t,d}\mathbf{X}_{k}+\mathbf{B}_{k'}\right),$$
 where $1\le k'\le m'$, while the sum index $k$ runs as below,
 $$k=\lfloor g(k'-1)/m\rfloor m/g+1,\dots,\left(\lfloor g(k'-1)/m\rfloor+1\right) m/g.$$
 Here,
 \begin{itemize}
     \item $\mathbf{W}\in\mathbb{R}^{m'\times (m/g)\times s}$ is the weight tensor
     \item $\otimes_{t,d}$ is the cross-correlation operator with stride $t$ and dilation $d$. That is, for any $\mathbf{w}\in\mathbb{R}^{s}$ and $\mathbf{x}\in\mathbb{R}^{n}$ one has $$\mathbf{w}\otimes_{t,d}\mathbf{x}\in\mathbb{R}^{\left\lfloor\frac{n-d(s-1)-1}{t} + 1\right\rfloor},$$ where
     $$\left(\mathbf{w}\otimes_{t,d}\mathbf{x}\right)_{j}:=\sum_{i=1}^{s}w_{i}x_{(j-1)t+(i-1)d+1}.$$
     \item $\mathbf{B}\in\mathbb{R}^{m'\times \left\lfloor\frac{n-d(s-1)-1}{t} + 1\right\rfloor}$ is the bias term.
 \end{itemize}
\end{definition}
The default values for stride and dilation are $t=1$, $d=1$. For this reason, with little abuse of notation, one says that $\Phi$ has no stride and no dilation to intend that $t=1$, $d=1$. Similarly, we assume $g=1$ whenever the grouping number is not declared explicitly.

\begin{definition}
\label{def:cnnt}
Let $m,m's,t,d$ be positive integers and let $g$ be a common divisor of $m$ and $m'$. A 1D Transposed Convolutional layer with $m$ input channels, $m'$ output channels, grouping number $g$, kernel size $s$, stride $t$, dilation factor $d$ and activation function $\rho:\mathbb{R}\to\mathbb{R}$, is a map of the form
$$\Phi:\mathbb{R}^{m\times n}\to\mathbb{R}^{m'\times (n-1)t+d(s-1)+1}$$
whose action on a given input $\mathbf{X}\in\mathbb{R}^{m\times n}$ is defined as
$$\Phi(\mathbf{X})_{k'}=\rho\left(\sum_{k}\mathbf{W}_{k,k'}\otimes_{t,d}^{\top}\mathbf{X}_{k}+\mathbf{B}_{k'}\right),$$
where $1\le k'\le m'$, while the sum index $k$ runs as below,
$$k=\lfloor g(k'-1)/m\rfloor m/g+1,\dots,\left(\lfloor g(k'-1)/m\rfloor+1\right) m/g.$$
Here,
 \begin{itemize}
     \item $\mathbf{W}\in\mathbb{R}^{(m/g)\times m'\times s}$ is the weight tensor
     \item $\otimes_{t,d}^{\top}$ is the transposed cross-correlation operator with stride $t$ and dilation $d$. That is, for any $\mathbf{w}\in\mathbb{R}^{s}$ and $\mathbf{x}\in\mathbb{R}^{n}$ one has $$\mathbf{w}\otimes_{t,d}^{\top}\mathbf{x}\in\mathbb{R}^{(n-1)t-d(s-1)+ 1},$$ where
     $$\left(\mathbf{w}\otimes_{t,d}^{\top}\mathbf{x}\right)_{j}:=\sum_{i}w_{\left\lfloor\frac{(i-1)t+1-j}{d}\right\rfloor+1}x_{i},$$
     the sum index $i$ running as below,
     $$i = \left\lfloor\frac{j-1}{t}+1\right\rfloor,\dots,\left\lfloor\frac{(s-1)d+j-1}{t}+1\right\rfloor.$$
     \item $\mathbf{B}\in\mathbb{R}^{m'\times (n-1)t-d(s-1)+1}$ is the bias term.
 \end{itemize}
\end{definition}

\begin{definition}
Let $\rho:\mathbb{R}\to\mathbb{R}$. A dense layer with activation function $\rho$ is a map $\Phi:\mathbb{R}^{n}\to\mathbb{R}^{n'}$ of the form
$$\Phi(\x)=\rho\left(\mathbf{W}\x+\mathbf{b}\right)$$
where $\mathbf{W}\in\mathbb{R}^{n'\times n}$ and $\mathbf{b}\in\mathbb{R}^{n'}$ are respectively the weight matrix and the bias vector.
\end{definition}

\begin{definition}
A Convolutional Neural Network (CNN) is any map that, up to reshaping operations, can be written as the composition of (transposed) convolutional layers. Conversely, a Deep Neural Network (DNN) is any map that, up to reshaping operations, can be written as the composition of dense layers.
\end{definition}
Since (transposed) convolutional layers can be seen as a particular class of dense layers, every CNN is a DNN. Similarly, CNNs and DNNs can be easily composed to build more complex DNN models.

\bibliographystyle{model1-num-names}
\bibliography{biblio}
\end{document}